\documentclass{article}
\usepackage{epsfig,epsfig}
\usepackage{amsmath}
\usepackage{amssymb}
\usepackage{booktabs,arydshln,multirow}
\usepackage{setspace}
\usepackage{xcolor}
\newtheorem{thm}{Theorem}[section]

\newtheorem{lemma}[thm]{Lemma}

\newtheorem{rem}[thm]{Remark}

\numberwithin{equation}{section} \topmargin=-1.5cm \oddsidemargin=-0.2cm
\begin{document}
\title{\textbf{TGMFE Algorithm Combined with Some Time Second-Order Schemes for Nonlinear Fourth-Order Reaction Diffusion System}
\thanks{\newline *Corresponding author. E-mail: mathliuyang@imu.edu.cn; mathliuyang@aliyun.com(Y. Liu)\newline \textit{Preprint submitted to Journal~~~~~~~~~~~~~~~~~~~~~~~~~~~~~~~~~~~~~~~~~~~~~~~~~~~~~~~~~~~~~~~~~~~~~~~~~~~~~~~~~~~~~~~~~~~~~~~~~~~~~~January 22, 2018}}}
\date{ }
\author{Baoli Yin$^1$, Yang Liu$^{1*}$, Hong Li$^1$, Siriguleng He$^1$, Jinfeng Wang$^2$
\\\small{1.School of Mathematical Sciences, Inner Mongolia University, Hohhot, 010021, China}\\\small {2. School of Statistics and Mathematics, Inner Mongolia University}\\\small { of Finance and Economics, Hohhot, 010070, China} }
\date{ }
 \maketitle
 {\color{black}\noindent\rule[0.25\baselineskip]{\textwidth}{0.25pt} }
\noindent\textbf{Abstract:}
%\begin{abstract}
%\noindent
In this article, a two-grid mixed finite element (TGMFE) method with some second-order time discrete schemes is developed for numerically solving
nonlinear fourth-order reaction diffusion equation. The two-grid MFE method is used to approximate spatial direction, and some second-order $\theta$ schemes formulated at time $t_{k-\theta}$ are considered to discretize the time direction. TGMFE method covers two main steps: a nonlinear MFE system based on the space coarse grid is solved by the iterative algorithm and a coarse solution is arrived at, then a linearized MFE system with fine grid is considered and a TGMFE solution is obtained. Here, the stability and a priori error estimates in $L^2$-norm for both nonlinear Galerkin MFE system and TGMFE scheme are derived. Finally, some convergence results are computed for both nonlinear Galerkin MFE system and TGMFE scheme to verify our theoretical analysis, which show that the convergence rate of the time second-order $\theta$ scheme including Crank-Nicolson scheme and second-order backward difference scheme is close to $2$, and that with the comparison to the computing time of nonlinear Galerkin MFE method, the CPU-time by using TGMFE method can be saved.
%\end{abstract}

{\color{black}\noindent\rule[0.25\baselineskip]{\textwidth}{0.25pt} }
%\\ \textbf{2000 Mathematics Subject Classification:} 65N30, 65M60, 26A33
\def\REF#1{\par\hangindent\parindent\indent\llap{#1\enspace}\ignorespaces}
\newcommand{\h}{\hspace{1.cm}}
\newcommand{\hh}{\hspace{2.cm}}
\newtheorem{yl}{\hspace{0.cm}Lemma}
\newtheorem{dl}{\hspace{0.cm}Theorem}
\newtheorem{re}{\hspace{0.cm}Remark}
\renewcommand{\sec}{\section*}
\renewcommand{\l}{\langle}
\renewcommand{\r}{\rangle}
\newcommand{\be}{\begin{eqnarray}}
\newcommand{\ee}{\end{eqnarray}}
\normalsize \vskip 0.2in

\section{Introduction}
Fourth-order partial differential equations (PDEs) have many important applications, so increasing people have made a lot of researches in looking for numerical solutions for these fourth-order PDEs. Li \cite{Li4} studied mixed methods with radial basis function for fourth-order parabolic and elliptic PDEs. In \cite{Li3}, Li gave the optimal convergence analysis of MFE methods for parabolic and elliptic problems with fourth-order derivatives. In \cite{Li2}, Li considered optimal error estimates of MFE method for a fourth-order elliptic problem with nonlinear term.
In \cite{Liuylh3}, Liu et al. discussed a MFE method for a fourth-order PDE with time-fractional derivative. In \cite{Khiari}, Khiari and Omrani considered a finite difference method for two-dimensional extended Fisher-Kolmogorov equation.
Liu et al. \cite{Liuylh2} considered MFE method for time-fractional fourth-order reaction-diffusion problem with a nonlinear term.
In \cite{C2}, Danumjaya and Pani considered $C^1$-conforming FE method for the extended Fisher-Kolmogorov (EFK) equation.
In \cite{Pani}, Danumjaya and Pani gave the studies of MFE methods for a fourth order reaction diffusion equation.
In \cite{XUDA1}, Yang et al. numerically solved fourth-order partial integro-differential equation by Crank-Nicolson/quasi-wavelets method.
Doss and Nandini \cite{DN} give a MFE method for
the extended Fisher-Kolmogorov equation. In \cite{Mohanty}, Mohanty and Kaur solved a class of fourth order non-linear parabolic equations by high accuracy implicit variable mesh methods. In \cite{Wangjf1}, Wang et al. solved the extended Fisher-Kolmogorov equation by using a new linearized Crank-Nicolson MFE scheme.
\par
In this article, we consider MFE method with two-grid scheme \cite{Xuj1,Xuj2} to look for the numerical solution of nonlinear fourth-order reaction diffusion equation
\begin{equation}\label{0.1}
 u_{t}+\gamma\bigtriangleup^2u-\bigtriangleup u+f(u)=g(\textbf{z},t), (\textbf{z},t)\in\Omega\times J,
\end{equation}
with initial condition
\begin{equation}\label{0.1c}
u(\textbf{z},0)=u_0(\textbf{z}), \textbf{z}\in\overline{\Omega},
\end{equation}
and  Dirichlet boundary conditions
\begin{equation}\label{0.1d}
u(\textbf{z},t)=\bigtriangleup u(\textbf{z},t)=0, (\textbf{z},t)\in\partial\Omega\times \bar{J},
\end{equation}
where $J=(0,T]$ is the time interval with $0<T<\infty$ and $\Omega$ is a bounded convex polygonal domain in $ R^d,d\leq 2$ with boundary $\partial\Omega$.
The coefficient $\gamma$ is a positive constant, $u_0(\textbf{z})$ is given function, $f(u)$ is the nonlinear term, which satisfies $|f(u_1)-f(u_2)|\leq C|u_1-u_2|$, $g(\textbf{z},t)$ is the source term. If we choose $\gamma=0$, the equation (\ref{0.1}) is the general second-order reaction diffusion equation; if we take $g(\textbf{z},t)=0$ and $f(u)=u^3-u$, the equation (\ref{0.1}) is called the extended Fisher-Kolmogorov (EFK) equation \cite{C2,Wangjf1}.
\par
Two-grid finite element (FE) method was proposed by Xu \cite{Xuj1,Xuj2}, which is an important numerical method and can save the computing time. This method includes mainly two computing steps: Firstly, a nonlinear system on the coarse grid is solved by the iterative method; secondly, based on the coarse solution solved by the first step, a linearized numerical scheme on the fine grid is constructed, then the two-grid solution for this linear system is arrived at. Compared to the nonlinear Galerkin FE method, the two-grid method can save the CPU time and also get the almost the same errors and convergence rate to the one of nonlinear Galerkin FE method. In view of the advantages of two-grid method, the method has been developed by increasing researchers, the detailed contents can be found in Dawson and Wheeler \cite{Dawson}, Chien and Jeng \cite{Chien}, Mu and Xu \cite{Xuj4}, Wu and Allen \cite{Wul1},
 Chen et al. \cite{Cheny1}, Chen and Chen \cite{Cheny3}, Liu et al. \cite{Ruih3}, Chen and Liu \cite{Ruih2}, Shi and Yang \cite{Shidy1}, Weng et al. \cite{Fengxl1}, Bajpai and Nataraj \cite{Baj}, Zhong et al. \cite{Zhonglq1}, Liu et al. \cite{Liuylh1},
 Liu et al. \cite{Liuydyw4}, Yan et al. \cite{Zhangzy1} and some other references. Based on these discussions for two-grid method, ones can see that the time direction is approximated mainly by the second-order Crank-Nicolson (CN) scheme, the second-order two step backward difference (BD) method, and backward Euler (BE) method with first-order convergence rate.

\par
Very recently, Liu et al. in \cite{Liuydu} developed some linearized time second-order $\theta$ schemes for nonlinear time fractional Cable equation. In \cite{Liuydu}, ones approximate the integer derivative $\frac{\partial u}{\partial t}(t_{n-\theta})$, $\theta\in [0,\frac12]$ by some second-order $\theta$ approximations.
Here we will develop the second-order $\theta$-scheme with TGMFE method for nonlinear fourth-order reaction diffusion equation, and we do not use the time linearized method for nonlinear term and only apply the second-order $\theta$-scheme to approximating time direction. In these $\theta$ schemes, our methods include general second-order CN method with $\theta=\frac12$ and second-order BD method with $\theta=0$.
For formulating our scheme, we introduce $\sigma=\bigtriangleup u$ and split the fourth-order problem (\ref{0.1}) into the coupled system of second-order equations
\begin{equation}\label{1.1}
 u_{t}+\gamma\bigtriangleup \sigma-\bigtriangleup u+f(u)=g(\textbf{z},t),
\end{equation}
and
\begin{equation}\label{1.2}
\sigma-\bigtriangleup u=0.
\end{equation}
In this paper, based on the lower coupled system, we consider the nonlinear Galerkin MFE system with some second-order $\theta$-scheme in time and second-order $\theta$ scheme with TGMFE method, then give the stability and error analysis with second-order convergence rate. For testing the effectiveness of second-order $\theta$ scheme and making some comparisons between nonlinear MFE system and TGMFE method, we provide some numerical examples.
\par
Throughout this article, $C>0$ is a constant, which is free of space mesh parameters $h$ and $H$ and the time step parameter $\Delta t$.
The layout of the paper is as follows. In Section 2, nonlinear Galerkin MFE scheme and TGMFE system are given. In Section 3, the stability and a priori error estimates are made. In Section 4, some numerical results are shown to confirm the theoretical analysis. In Section 5, some conclusions are shown.
\section{Numerical scheme}
For obtaining fully discrete scheme, we insert the nodes $t_n=n\Delta t (n=0,1,2,\cdots,N)$ in the time
interval $[0,T]$, where $t_n$ satisfy $0=t_0<t_1<t_2<\cdots<t_N=T$ with mesh length $\Delta t=T/N$
for some positive integer $N$. We now define $\phi^n$ as the approximation of $\phi(t_n)$ at time $t=t_n$.
\par
To formulate the time discrete scheme, referring to Ref. \cite{Liuydu}, we have the following lemmas at time $t=t_{n-\theta}$.
\begin{lemma}\label{lemma1}For sufficiently smooth function $\phi(t)$, at time $t_{n-\theta}$, the following approximation for first-order derivative with second-order convergence rate for any $\theta\in [0,\frac12]$ holds
\begin{equation}\begin{split}\label{4.6}
\phi_t(t_{n-\theta})=&\frac{(3-2\theta)\phi(t_{n})-(4-4\theta)\phi(t_{n-1})+(1-2\theta)\phi(t_{n-2})}{2\Delta t}+O(\Delta t^2)
\\=&\frac{(3-2\theta)\phi^{n}-(4-4\theta)\phi^{n-1}+(1-2\theta)\phi^{n-2}}{2\Delta t}\\\triangleq&\mathcal{D}_t\phi^{n-\theta},~n\geq2.
\end{split}\end{equation}
\end{lemma}
\begin{lemma}\label{lemma5}With sufficiently smooth function $\phi(t)$, at time $t_{n-\theta}$, the following approximate formula
\begin{equation}\begin{split}\label{4.3}
\phi(t_{n-\theta})=&(1-\theta)\phi(t_{n})+\theta \phi(t_{n-1})+O(\Delta t^2)
\\
=&(1-\theta)\phi^{n}+\theta \phi^{n-1}
\\\triangleq& \phi^{n-\theta}
\end{split}\end{equation}
holds for any $\theta\in [0,\frac12]$.
\end{lemma}
\begin{lemma}\label{lemma6}For series $\{\phi^n\}$ ($n\geq2$), the following inequality holds
\begin{equation}\label{4.15}\Big{(}\mathcal{D}_t\phi^{n-\theta},\phi^{n-\theta}\Big{)}\geq \frac{1}{4\Delta t}(\mathbb{H}[\phi^{n}]-\mathbb{H}[\phi^{n-1}])
,
\end{equation}
\begin{equation}\label{4.14}\begin{split}
&\mathbb{H}[\phi^{n}]=(3-2\theta){\|\phi^{n}\|}^{2}-(1-2\theta){\|\phi^{n-1}\|}^{2}+(2-\theta)(1-2\theta){\|\phi^{n}-\phi^{n-1}\|}^{2},~n\geq1,
\end{split}
\end{equation}
\text{and}
\begin{equation}\label{3.2}\mathbb{H}[\phi^{n}]\geq\frac{1}{1-\theta}\|\phi^{n}\|^{2},\end{equation}
where $0\leq\theta\leq1/2$.
\end{lemma}
Based on the above lemmas \ref{lemma1}-\ref{lemma5}, we take the values for system (\ref{1.1})-(\ref{1.2}) at time $t_{n-\theta}$ and drop the time errors  to get for $n\geq2$
\begin{equation}\label{1.12}
\mathcal{D}_tu^{n-\theta}+\gamma\bigtriangleup \sigma^{n-\theta}-\bigtriangleup u^{n-\theta}+f^{n-\theta}(u)=g(\textbf{z},t_{n-\theta}),
\end{equation}
and
\begin{equation}\label{1.22}
\sigma^{n-\theta}-\bigtriangleup u^{n-\theta}=0.
\end{equation}
Similarly, we have for $n=1$
\begin{equation}\label{1.21}
\mathcal{D}_tu^{\frac12}+\gamma\bigtriangleup \sigma^{\frac12}-\bigtriangleup u^{\frac12}+f^{\frac12}(u)=g(\textbf{z},t_{\frac12}),
\end{equation}
and
\begin{equation}\label{1.23}
\sigma^{\frac12}-\bigtriangleup u^{\frac12}=0.
\end{equation}
At time $t_{n-\theta}$, we use lemmas \ref{lemma1}-\ref{lemma5} to get time semi-discrete scheme for
the case $n\geq2$
\begin{equation}\begin{split}\label{2.2}
\Big{(}\mathcal{D}_tu^{n-\theta},v\Big{)}-\gamma(\nabla \sigma^{n-\theta},\nabla v)+(\nabla u^{n-\theta},\nabla v)+(f^{n-\theta}(u),v)=&(g^{n-\theta},v),~\forall v\in H_0^1,
\end{split}
\end{equation}
and
\begin{equation}\label{2.3}
(\sigma^{n-\theta},w)+(\nabla u^{n-\theta},\nabla w)=0,~\forall w\in H_0^1.
\end{equation}
For solving time semi-discrete coupled system (\ref{2.2})-(\ref{2.3}), we need to solve the solution pair $(u^1,\sigma^1)$ based on the initial value $(u^0,\sigma^0)$ by the following time discrete system
\begin{equation}\begin{split}\label{2.21}
\Big{(}\mathcal{D}_tu^{\frac12},v\Big{)}-\gamma(\nabla \sigma^{\frac12},\nabla v)+(\nabla u^{\frac12},\nabla v)+(f^{\frac12}(u),v)=&(g^{\frac12},v),~\forall v\in H_0^1,
\end{split}
\end{equation}
and
\begin{equation}\label{2.31}
(\sigma^{\frac12},w)+(\nabla u^{\frac12},\nabla w)=0,~\forall w\in H_0^1.
\end{equation}
Considering the above time semi-discrete coupled system (\ref{2.2})-(\ref{2.31}), we get standard nonlinear MFE scheme for
the case $n\geq2$
\begin{equation}\begin{split}\label{2.4}
\Big{(}\mathcal{D}_tu_h^{n-\theta},v_h\Big{)}-\gamma(\nabla \sigma_h^{n-\theta},\nabla v_h)+(\nabla u_h^{n-\theta},\nabla v_h)+(f^{n-\theta}(u_h),v_h)=&(g^{n-\theta},v_h),~\forall v_h\in L_h,
\end{split}
\end{equation}
and
\begin{equation}\label{2.5}
(\sigma_h^{n-\theta},w_h)+(\nabla u_h^{n-\theta},\nabla w_h)=0,~\forall w_h\in L_h,
\end{equation}
in which the MFE solution pair $(u_h^1,\sigma^1_h)$ can be arrived at by the following coupled MFE scheme
\begin{equation}\begin{split}\label{2.41}
\Big{(}\mathcal{D}_tu_h^{\frac12},v_h\Big{)}-\gamma(\nabla \sigma_h^{\frac12},\nabla v_h)+(\nabla u_h^{\frac12},\nabla v_h)+(f^{\frac12}(u_h),v_h)=&(g^{\frac12},v_h),~\forall v_h\in L_h,
\end{split}
\end{equation}
and
\begin{equation}\label{2.51}
(\sigma_h^{\frac12},w_h)+(\nabla u_h^{\frac12},\nabla w_h)=0,~\forall w_h\in L_h.
\end{equation}
where $L_h$ is the finite element subspace of $H_0^1$.
\par
For saving the CPU-time of standard MFE method, we consider TGMFE method based on the coarse grid $\mathfrak{T}_H$ and the fine grid $\mathfrak{T}_h$ for $n\geq1$.
\\
\textbf{Step I}: Based on the coarse grid $\mathfrak{T}_H$, we solve nonlinear coupled system by looking for the coarse solution pair $\{u^{n}_H,\sigma^{n}_H\}:[0,T]\mapsto L_H\subset L_h$ such that for the case $n\geq2$:
\begin{equation}\begin{split}\label{2.13}
\Big{(}\mathcal{D}_tu_H^{n-\theta},v_H\Big{)}-\gamma(\nabla \sigma_H^{n-\theta},\nabla v_H)+(\nabla u_H^{n-\theta},\nabla v_H)+(f^{n-\theta}(u_H),v_H)=&(g^{n-\theta},v_H),~\forall v_H\in L_H,
\end{split}
\end{equation}
and
\begin{equation}\label{2.14}
(\sigma_H^{n-\theta},w_h)+(\nabla u_H^{n-\theta},\nabla w_H)=0,~\forall w_H\in L_H,
\end{equation}
which can be solved based on the coarse solution pairs $(u_H^1,\sigma^1_H)$ and $(u_H^0,\sigma^0_H)$, where $(u_H^1,\sigma^1_H)$ need to be solved by the following system
\begin{equation}\begin{split}\label{2.131}
\Big{(}\mathcal{D}_tu_H^{\frac12},v_H\Big{)}-\gamma(\nabla \sigma_H^{\frac12},\nabla v_H)+(\nabla u_H^{\frac12},\nabla v_H)+(f^{\frac12}(u_H),v_H)=&(g^{\frac12},v_H),~\forall v_H\in L_H,
\end{split}
\end{equation}
and
\begin{equation}\label{2.141}
(\sigma_H^{\frac12},w_h)+(\nabla u_H^{\frac12},\nabla w_H)=0,~\forall w_H\in L_H.
\end{equation}
\textbf{Step II}: Based on the coarse solution pair $\{u^{n}_H,\sigma^{n}_H\}\in L_H$, we solve a linear coupled system on the fine grid $\mathfrak{T}_h$, by finding the fine solution pair $\{U^{n}_h,\Sigma^{n}_h\}:[0,T]\mapsto L_h$ such that for the case $n\geq2$:
\begin{equation}\begin{split}\label{2.15}
\Big{(}\mathcal{D}_tU_h^{n-\theta},v_h\Big{)}&-\gamma(\nabla \Sigma_h^{n-\theta},\nabla v_h)+(\nabla U_h^{n-\theta},\nabla v_h)\\&+((1-\theta)\mathfrak{F}(U_{h}^{n},u_{H}^{n})+\theta f(U_{h}^{n-1}),v_h)=(g^{n-\theta},v_h),~\forall v_h\in L_h,
\end{split}
\end{equation}
and
\begin{equation}\label{2.16}
(\Sigma_h^{n-\theta},w_h)+(\nabla U_h^{n-\theta},\nabla w_h)=0,~\forall w_h\in L_h,
\end{equation}
in which the two-grid solution pair $(U_h^1,\Sigma^1_h)$ can be arrived at by the system
\begin{equation}\begin{split}\label{2.151}
\Big{(}\mathcal{D}_tU_h^{\frac12},v_h\Big{)}&-\gamma(\nabla \Sigma_h^{\frac12},\nabla v_h)+(\nabla U_h^{\frac12},\nabla v_h)\\&+\frac12(\mathfrak{F}(U_{h}^{1},u_{H}^{1})+ f(U_{h}^{0}),v_h)=(g^{\frac12},v_h),~\forall v_h\in L_h,
\end{split}
\end{equation}
and
\begin{equation}\label{2.161}
(\Sigma_h^{\frac12},w_h)+(\nabla U_h^{\frac12},\nabla w_h)=0,~\forall w_h\in L_h,
\end{equation}
where $\mathfrak{F}(U_{h}^{n},u_{H}^{n})\triangleq f(u_{H}^{n})+f'(u_{H}^{n})(U_{h}^{n}-u_{H}^{n})$, $n\geq1$.
\begin{rem}(i).
With $\sigma=\Delta u$, we have the second mixed formulation
for the case $n\geq2$
\begin{equation}\begin{split}\label{2.111}
\Big{(}\mathcal{D}_tu^{n-\theta},v\Big{)}-\gamma(\nabla \sigma^{n-\theta},\nabla v)+(\sigma^{n-\theta},v)+(f^{n-\theta}(u),v)=&(g^{n-\theta},v),~\forall v\in H_0^1,
\end{split}
\end{equation}
\begin{equation}\label{2.112}
(\sigma^{n-\theta},w)+(\nabla u^{n-\theta},\nabla w)=0,~\forall w\in H_0^1.
\end{equation}
(ii).
With $\sigma=\gamma\Delta u-u$, we have the third mixed formulation
for the case $n\geq2$
\begin{equation}\begin{split}\label{2.113}
\Big{(}\mathcal{D}_tu^{n-\theta},v\Big{)}-(\nabla \sigma^{n-\theta},\nabla v)+(f^{n-\theta}(u),v)=&(g^{n-\theta},v),~\forall v\in H_0^1,
\end{split}
\end{equation}
\begin{equation}\label{2.114}
(\sigma^{n-\theta},w)+\gamma(\nabla u^{n-\theta},\nabla w)+(u^{n-\theta},w)=0,~\forall w\in H_0^1.
\end{equation}
(iii). For the above two schemes (i)-(ii), we can discuss the numerical theories by the similar process to the one of the current scheme. Here, we do not consider that.
\end{rem}
\section{Numerical analysis for fully discrete scheme}
\subsection{Stability}
In what follows, we consider the following stable inequality.
\begin{thm}\label{thm1}For the standard MFE solution pair $\{u_h^n,\sigma_h^n\}\in L_h$, the stability for nonlinear MFE system (\ref{2.4})-(\ref{2.5}) holds
\begin{equation}\begin{split}\label{2.20}
\|u_h^{n}\|^2+\Delta t\sum_{k=2}^n\|\nabla u_h^{k-\theta}\|^2+\gamma\Delta t\sum_{k=2}^n\|\sigma_h^{k-\theta}\|^2
\leq& C(\|u_h^{0}\|^2+\Delta t\sum_{k=1}^n\|g^{k}\|^2).
\end{split}
\end{equation}
\end{thm}
\textbf{Proof}. In (\ref{2.4})-(\ref{2.5}), we choose $v_{h}=u_h^{n-\theta}$ and $w_{h}=\sigma_h^{n-\theta}$ to arrive at
\begin{equation}\begin{split}\label{2.8}
\Big{(}\mathcal{D}_tu_h^{n-\theta},u_h^{n-\theta}\Big{)}+\|\nabla u_h^{n-\theta}\|^2=&-\gamma(\nabla \sigma_h^{n-\theta},\nabla u_h^{n-\theta})-(f^{n-\theta}(u_h),u_h^{n-\theta})+(g^{n-\theta},v_h)
\\=&-\gamma\|\sigma_h^{n-\theta}\|^2-(f^{n-\theta}(u_h),u_h^{n-\theta})+(g^{n-\theta},u_h^{n-\theta}).
\end{split}
\end{equation}
We use Cauchy-Schwarz inequality as well as Young inequality to get
\begin{equation}\begin{split}\label{2.9}
\frac{1}{4\Delta t}(\mathbb{H}[u_h^{n}]-\mathbb{H}[u_h^{n-1}])+\|\nabla u_h^{n-\theta}\|^2+\gamma\|\sigma_h^{n-\theta}\|^2
\leq&-(f^{n-\theta}(u_h),u_h^{n-\theta})+(g^{n-\theta},u_h^{n-\theta})
\\\leq&\frac12\|f^{n-\theta}(u_h)\|^2+\|u_h^{n-\theta}\|^2+\frac12\|g^{n-\theta}\|^2.
\end{split}
\end{equation}
Sum from $2$ to $n$ for the above inequality and make use of lemma \ref{lemma6} to get
\begin{equation}\begin{split}\label{2.10}
&\mathbb{H}[u_h^{n}]+4\Delta t\sum_{k=2}^n\|\nabla u_h^{k-\theta}\|^2+4\gamma\Delta t\sum_{k=2}^n\|\sigma_h^{k-\theta}\|^2
\\\leq&\mathbb{H}[u_h^{2}]+2\Delta t\sum_{k=2}^n\|(1-\theta)u_h^{k}+\theta u_h^{k-1}\|^2+2\Delta t\sum_{k=2}^n\|(1-\theta)g^{k}+\theta g^{k-1}\|^2
\\\leq&\mathbb{H}[u_h^{2}]+2\Delta t\sum_{k=1}^n\|u_h^{k}\|^2+2\Delta t\sum_{k=1}^n\|g^{k}\|^2.
\end{split}
\end{equation}
In the next step, we need to estimate $\mathbb{H}[u_h^{2}]$. We first compute $u_h^1$ by using Crank-Nicolson scheme, then apply lemma \ref{lemma6}, we have
\begin{equation}\begin{split}\label{2.11}
&\frac{1}{1-\theta}\|u_h^{n}\|^2+4\Delta t\sum_{k=2}^n\|\nabla u_h^{k-\theta}\|^2+4\gamma\Delta t\sum_{k=2}^n\|\sigma_h^{k-\theta}\|^2
\\\leq&C\Delta t\sum_{k=0}^n\|u_h^{k}\|^2+2\Delta t\sum_{k=1}^n\|g^{k}\|^2.
\end{split}
\end{equation}
Using Gronwall lemma for the above inequality, we complete the proof.
\begin{thm}\label{thm2}For the coarse solution pair $\{u_H^n,\sigma_H^n\}\in L_H$, the stability for the coupled system (\ref{2.13})-(\ref{2.14}) holds
\begin{equation}\begin{split}\label{2.12}
\|u_H^{n}\|^2+\Delta t\sum_{k=2}^n\|\nabla u_H^{k-\theta}\|^2+\gamma\Delta t\sum_{k=2}^n\|\sigma_H^{k-\theta}\|^2
\leq& C(\|u_H^{0}\|^2+\Delta t\sum_{k=1}^n\|g^{k}\|^2).
\end{split}
\end{equation}
\end{thm}
\textbf{Proof.} Follow the similar process to the one in theorem \ref{thm1} to accomplish the proof.
\par
In what follows, we will give the stability of TGMFE scheme.
\begin{thm}\label{thm3}For the two-grid solution pair $\{U_h^n,\Sigma_h^n\}\in L_h$, the stability for the TGMFE system (\ref{2.13})-(\ref{2.161}) holds
\begin{equation}\begin{split}\label{2.221}
\|U_h^{n}\|^2+\Delta t\sum_{k=2}^n\|\nabla U_h^{k-\theta}\|^2+\gamma\Delta t\sum_{k=2}^n\|\Sigma_h^{k-\theta}\|^2
\leq& C(\|U_h^{0}\|^2+\|u_H^{0}\|^2+\Delta t\sum_{k=1}^n\|g^{k}\|^2).
\end{split}
\end{equation}
\end{thm}
\textbf{Proof.} Using the similar derivation to that in theorem \ref{thm1} again, we have
\begin{equation}\begin{split}\label{2.23}
&\|U_h^{n}\|^2+\Delta t\sum_{k=2}^n\|\nabla U_h^{k-\theta}\|^2+\gamma\Delta t\sum_{k=2}^n\|\Sigma_h^{k-\theta}\|^2
\\\leq& C(\|U_h^{0}\|^2+\Delta t\sum_{k=1}^n\|g^{k}\|^2)+\Delta t\sum_{k=1}^n\|U_h^{n}\|^2+\Delta t\sum_{k=1}^n\|u_H^{n}\|^2.
\end{split}
\end{equation}
Combining Gronwall lemma with (\ref{2.12}), we obtain the result (\ref{2.221}).
\subsection{Error analysis}
To carry out the error analysis, we give the projection operator and the estimate inequality.
\begin{lemma}\label{lemma7}
Define a Ritz projection operator $\mathfrak{Q}_{\mathfrak{h}}:H^{1}_{0}(\Omega)\rightarrow L_{\mathfrak{h}}$ satisfying\\
\begin{equation}\label{2.26}
(\nabla (y-\mathfrak{Q}_{\mathfrak{h}}y),\nabla y_{\mathfrak{h}})=0,~\forall y_{\mathfrak{h}}\in V_{\mathfrak{h}},
\end{equation}
with the estimate inequality
\begin{equation}\label{2.25}
\|y-\mathfrak{Q}_{\mathfrak{h}}y\|+\|(y-\mathfrak{Q}_{\mathfrak{h}}y)_t\|+\mathfrak{h}\|y-\mathfrak{Q}_{\mathfrak{h}}y\|_{1}\leq C\mathfrak{h}^{m+1},~\forall y\in H^1_0(\Omega)\cap H^{m+1}(\Omega),
\end{equation}
where the norms are defined by $\|y\|_{l}=\sqrt{\sum_{0\leq\mid r\mid\leq l}\int_{\Omega}\mid D^{r}y\mid^2}$ and $\mathfrak{h}=h~\text{or}~H$.
\end{lemma}
In what follows, we will give the detailed proof of error estimates in $L^2$-norm.
%\begin{thm}\label{thm4}Let $u^n$ be the solution of system (\ref{2.2})-(\ref{2.31}) and $u_h^n$ be the solution of system (\ref{2.4})-(\ref{2.51}), %respectively. With $u_{h}^0=\mathfrak{Q}_{h}u_0$, there exists a constant $C$ free of space-time mesh step pair $(h,\Delta t)$ such that
%\begin{equation}\label{5.16}
%\|u^n-u_h^{n}\|\leq Ch^{m+1}.
%\end{equation}
%\end{thm}
%\textbf{Proof}.
\begin{thm}\label{thm5}Let $u^n$ be the solution of system (\ref{2.2})-(\ref{2.31}), $u_h^n$ be the solution of system (\ref{2.4})-(\ref{2.51}), $u_H^n$ be the coarse solution of system (\ref{2.13})-(\ref{2.141}), $U_h^n$ be the two-grid solution of system (\ref{2.13})-(\ref{2.161}), respectively. With $u_{h}^0=\mathfrak{Q}_{h}u_0$, $u_{H}^0=\mathfrak{Q}_{H}u_0$, $U_{h}^0=\mathfrak{Q}_{h}u_0$, there exists a constant $C$ free of space-time mesh step length $h,H,\Delta t$ such that
\begin{equation}\label{5.16}
\|u^n-u_h^{n}\|+\gamma^{\frac12}\Big{(}\Delta t\sum_{k=1}^n\|(\sigma-\sigma_h)^{k-\theta}\|^2\Big{)}^{\frac12}+\Big{(}\Delta t\sum_{k=1}^n\|\nabla (u-u_{h})^{k-\theta}\|^2\Big{)}^{\frac12}\leq Ch^{m+1},
\end{equation}
\begin{equation}\label{5.17}
\|u^n-u_H^{n}\|+\gamma^{\frac12}\Big{(}\Delta t\sum_{k=1}^n\|(\sigma-\sigma_H)^{k-\theta}\|^2\Big{)}^{\frac12}+\Big{(}\Delta t\sum_{k=1}^n\|\nabla (u-u_{H})^{k-\theta}\|^2\Big{)}^{\frac12}\leq CH^{m+1},
\end{equation}
and
\begin{equation}\label{5.18}
\|u^n-U_h^{n}\|+\gamma^{\frac12}\Big{(}\Delta t\sum_{k=1}^n\|(\sigma-\Sigma_h)^{k-\theta}\|^2\Big{)}^{\frac12}+\Big{(}\Delta t\sum_{k=1}^n\|\nabla (u-U_{h})^{k-\theta}\|^2\Big{)}^{\frac12}\leq C[h^{m+1}+H^{2m+2}].
\end{equation}
\end{thm}
\textbf{Proof.} \textbf{1)}. Combine (\ref{2.2})-(\ref{2.3}) with (\ref{2.4})-(\ref{2.5}), we use projection operator (\ref{2.26}) to get
\begin{equation}\begin{split}\label{2.215}
\Big{(}\mathcal{D}_t(u-u_h)^{n-\theta},v_h\Big{)}&-\gamma(\nabla (\mathfrak{Q}_h\sigma-\sigma_h)^{n-\theta},\nabla v_h)\\&+(\nabla (\mathfrak{Q}_hu-u_h)^{n-\theta},\nabla v_h)+(f^{n-\theta}(u)-f^{n-\theta}(u_h),v_h)=0,
\end{split}
\end{equation}
and
\begin{equation}\label{2.22}
((\sigma-\mathfrak{Q}_h\sigma)^{n-\theta},w_h)+(\nabla (\mathfrak{Q}_hu-u_h)^{n-\theta},\nabla w_h)=0.
\end{equation}
Now we take $(v_h,w_h)=((\mathfrak{Q}_hu-u_h)^{n-\theta},(\mathfrak{Q}_h\sigma-\sigma_h)^{n-\theta})$ in system (\ref{2.215})-(\ref{2.22}) to get
\begin{equation}\begin{split}\label{2.23}
&\Big{(}\mathcal{D}_t(\mathfrak{Q}_hu-u_h)^{n-\theta},(\mathfrak{Q}_hu-u_h)^{n-\theta}\Big{)}+\gamma\|(\mathfrak{Q}_h\sigma-\sigma_h)^{n-\theta}\|^2
+\|\nabla (\mathfrak{Q}_hu-u_h)^{n-\theta}\|^2\\
=&-(f^{n-\theta}(u)-f^{n-\theta}(u_h),(\mathfrak{Q}_hu-u_h)^{n-\theta})-\Big{(}\mathcal{D}_t(u-\mathfrak{Q}_hu)^{n-\theta},(\mathfrak{Q}_hu-u_h)^{n-\theta}\Big{)}
\\&-\gamma((\sigma-\mathfrak{Q}_h\sigma)^{n-\theta},(\mathfrak{Q}_h\sigma-\sigma_h)^{n-\theta}).
\end{split}
\end{equation}
Now we estimate the three terms on the right hand side of the above equation. For the first term, we use lemma \ref{lemma6}, Cauchy-Schwarz inequality as well as Young inequality to get
\begin{equation}\begin{split}\label{2.24}
&-(f^{n-\theta}(u)-f^{n-\theta}(u_h),(\mathfrak{Q}_hu-u_h)^{n-\theta})
\\\leq& \|f^{n-\theta}(u)-f^{n-\theta}(u_h)\|\|(\mathfrak{Q}_hu-u_h)^{n-\theta})\|
\\\leq& \|(1-\theta)(f^{n}(u)-f^{n}(u_h))+\theta (f^{n-1}(u)-f^{n-1}(u_h))\|\|(\mathfrak{Q}_hu-u_h)^{n-\theta})\|
\\\leq& C(\|u^n-\mathfrak{Q}_hu^n\|^2+\|\mathfrak{Q}_hu^n-u_h^n\|^2+ \|u^{n-1}-\mathfrak{Q}_hu^{n-1}\|^2+\|\mathfrak{Q}_hu^{n-1}-u_h^{n-1}\|^2).
\end{split}
\end{equation}
For the second and third terms on the right hand side of (\ref{2.23}), we use Cauchy-Schwarz inequality as well as Young inequality to arrive at
\begin{equation}\begin{split}\label{2.48}
&-\Big{(}\mathcal{D}_t(u-\mathfrak{Q}_hu)^{n-\theta},(\mathfrak{Q}_hu-u_{h})^{n-\theta}\Big{)}
-\gamma((\sigma-\mathfrak{Q}_h\sigma)^{n-\theta},(\mathfrak{Q}_h\sigma-\sigma_{h}))^{n-\theta})\\\leq&
\frac12\|\mathcal{D}_t(u-\mathfrak{Q}_hu)^{n-\theta}\|^2+\frac12\|(\mathfrak{Q}_hu-u_{h})^{n-\theta}\|^2
+\frac{\gamma}{2}\|(\sigma-\mathfrak{Q}_h\sigma)^{n-\theta}\|^2+\frac{\gamma}{2}\|(\mathfrak{Q}_h\sigma-\sigma_{h})^{n-\theta}\|^2
\\\leq&\frac12\Big{\|}\frac{(3-2\theta)[(u-\mathfrak{Q}_hu)^{n}-(u-\mathfrak{Q}_hu)^{n-1}]
-(1-2\theta)[(u-\mathfrak{Q}_hu)^{n-1}-(u-\mathfrak{Q}_hu)^{n-2}]}{2\Delta t}\Big{\|}^2\\&+\frac12\|(\mathfrak{Q}_hu-U_{h})^{n-\theta}\|^2
+\frac{\gamma}{2}\|(\sigma-\mathfrak{Q}_h\sigma)^{n-\theta}\|^2+\frac{\gamma}{2}\|(\mathfrak{Q}_h\sigma-\sigma_{h})^{n-\theta}\|^2
\\\leq&\frac{C}{\Delta t}\int_{t_{n-2}}^{t_n}\|(u-\mathfrak{Q}_hu)_t\|^2dt+\frac12\|(\mathfrak{Q}_hu-u_{h})^{n-\theta}\|^2
+\frac{\gamma}{2}\|(\sigma-\mathfrak{Q}_h\sigma)^{n-\theta}\|^2+\frac{\gamma}{2}\|(\mathfrak{Q}_h\sigma-\sigma_{h})^{n-\theta}\|^2.
\end{split}
\end{equation}
Substitute (\ref{2.24}) and (\ref{2.48}) into (\ref{2.23}) and use inequality (\ref{4.15}) to get
\begin{equation}\begin{split}\label{2.49}
&\frac{1}{4\Delta t}(\mathbb{H}[(\mathfrak{Q}_hu-u_{h})^{n}]-\mathbb{H}[(\mathfrak{Q}_hu-u_{h})^{n-1}])
+\frac{\gamma}{2}\|(\mathfrak{Q}_h\sigma-\sigma_h)^{n-\theta}\|^2+\|\nabla (\mathfrak{Q}_hu-u_{h})^{n-\theta}\|^2
\\\leq&\frac{C}{\Delta t}\int_{t_{n-2}}^{t_n}\|(u-\mathfrak{Q}_hu)_t\|^2dt+\frac12\|(\mathfrak{Q}_hu-u_{h})^{n-\theta}\|^2
+\frac{\gamma}{2}\|(\sigma-\mathfrak{Q}_h\sigma)^{n-\theta}\|^2
\\&+C(\|u^n-u_{h}^{n}\|^2+\|u^{n-1}-u_{h}^{n-1}\|^2+\|(\mathfrak{Q}_hu-u_{h})^{n-\theta}\|^2).
\end{split}
\end{equation}
Sum (\ref{2.49}) with respect to $n$ from $2$ to $L$ to get
\begin{equation}\begin{split}\label{2.50}
&\mathbb{H}[(\mathfrak{Q}_hu-u_{h})^{L}]
+2\Delta t\sum_{n=2}^L(\gamma\|(\mathfrak{Q}_h\sigma-\sigma_h)^{n-\theta}\|^2+2\|\nabla (\mathfrak{Q}_hu-u_{h})^{n-\theta}\|^2)
\\\leq&\mathbb{H}[(\mathfrak{Q}_hu-u_{h})^{1}]+C\int_{t_{0}}^{t_L}\|(u-\mathfrak{Q}_hu)_t\|^2dt+2\Delta t\sum_{n=2}^L\|(\mathfrak{Q}_hu-u_{h})^{n-\theta}\|^2
+2\gamma\Delta t\sum_{n=2}^L\|(\sigma-\mathfrak{Q}_h\sigma)^{n-\theta}\|^2
\\&+C\Delta t\sum_{n=2}^L(\|u^n-u_{h}^{n}\|^2
+\|u^{n-1}-u_{h}^{n-1}\|^2+\|(\mathfrak{Q}_hu-u_{h})^{n-\theta}\|^2).
\end{split}
\end{equation}
Now we need to give the estimate of $\mathbb{H}[(\mathfrak{Q}_hu-u_{h})^{1}]$. We combine (\ref{2.21})-(\ref{2.31}) with (\ref{2.41})-(\ref{2.51}) to get
\begin{equation}\begin{split}\label{2.52}
\Big{(}\mathcal{D}_t(u-u_h)^{\frac12},v_h\Big{)}&-\gamma(\nabla (\sigma-\sigma_h)^{\frac12},\nabla v_h)+(\nabla (u-u_h)^{\frac12},\nabla v_h)\\&+(f^{\frac12}(u)-f^{\frac12}(u_{h}),v_h)=0,
\end{split}
\end{equation}
and
\begin{equation}\label{2.53}
((\sigma-\sigma_h)^{\frac12},w_h)+(\nabla (u-u_h)^{\frac12},\nabla w_h)=0.
\end{equation}
We take $(v_h,w_h)=((\mathfrak{Q}_hu-u_{h})^{\frac12},(\mathfrak{Q}_h\sigma-\sigma_{h})^{\frac12})$ in (\ref{2.52})-(\ref{2.53}), and use the similar derivation to the one of inequality (\ref{2.50}) to get
\begin{equation}\begin{split}\label{2.54}
&\|(\mathfrak{Q}_hu-u_{h})^{1}\|^2
+2\Delta t(\gamma\|(\mathfrak{Q}_h\sigma-\sigma_h)^{\frac12}\|^2+2\|\nabla (\mathfrak{Q}_hu-u_{h})^{\frac12}\|^2)
\\\leq&\|(\mathfrak{Q}_hu-u_{h})^{0}\|^2+C\int_{t_{0}}^{t_1}\|(u-\mathfrak{Q}_hu)_t\|^2dt+2\Delta t\|(\mathfrak{Q}_hu-u_{h})^{\frac12}\|^2
+2\gamma\Delta t\|(\sigma-\mathfrak{Q}_h\sigma)^{\frac12}\|^2
\\&+C\Delta t(\|u^1-u_{h}^{1}\|^2
+\|u^{0}-u_{h}^{0}\|^2+\|(\mathfrak{Q}_hu-u_{h})^{\frac12}\|^2).
\end{split}
\end{equation}
From (\ref{2.54}), we easily know that
\begin{equation}\begin{split}\label{2.55}
\mathbb{H}[(\mathfrak{Q}_hu-u_{h})^{1}]=&(3-2\theta){\|(\mathfrak{Q}_hu-u_{h})^{1}\|}^{2}-(1-2\theta){\|(\mathfrak{Q}_hu-u_{h})^{0}\|}^{2}
\\&+(2-\theta)(1-2\theta){\|(\mathfrak{Q}_hu-u_{h})^{1}-(\mathfrak{Q}_hu-u_{h})^{0}\|}^{2}\\\leq & C(\|(\mathfrak{Q}_hu-u_{h})^{1}\|^2+\|(\mathfrak{Q}_hu-u_{h})^{0}\|^2).
\end{split}
\end{equation}
Combine (\ref{2.50}), (\ref{2.54}) with (\ref{2.55}) to get
\begin{equation}\begin{split}\label{2.56}
&\mathbb{H}[(\mathfrak{Q}_hu-u_{h})^{L}]
+2\Delta t\sum_{n=1}^L(\gamma\|(\mathfrak{Q}_h\sigma-\sigma_h)^{n-\theta}\|^2+2\|\nabla (\mathfrak{Q}_hu-u_{h})^{n-\theta}\|^2)
\\\leq&C\int_{t_{0}}^{t_L}\|(u-\mathfrak{Q}_hu)_t\|^2dt+2\Delta t\sum_{n=2}^L\|(\mathfrak{Q}_hu-u_{h})^{n-\theta}\|^2
+2\gamma\Delta t\sum_{n=2}^L\|(\sigma-\mathfrak{Q}_h\sigma)^{n-\theta}\|^2
\\&+C\Delta t\sum_{n=2}^L(\|u^n-u_{h}^{n}\|^2
+\|u^{n-1}-u_{h}^{n-1}\|^2+\|(\mathfrak{Q}_hu-u_{h})^{n-\theta}\|^2),
\end{split}
\end{equation}
which is combined with Gronwall lemma, triangle inequality and (\ref{2.25}) to arrive at the conclusion (\ref{5.16}).
\\
\textbf{2)}. Using the similar proof to the one of (\ref{5.16}), we arrive at the conclusion (\ref{5.17}) based on the coarse grid.
\\
\textbf{3)}. Now we combine system (\ref{2.2})-(\ref{2.3}) with (\ref{2.15})-(\ref{2.16}) to get for $n\geq 2$
\begin{equation}\begin{split}\label{2.33}
\Big{(}\mathcal{D}_t(u-U_h)^{n-\theta},v_h\Big{)}&-\gamma(\nabla (\sigma-\Sigma_h)^{n-\theta},\nabla v_h)+(\nabla (u-U_h)^{n-\theta},\nabla v_h)\\&+(f^{n-\theta}(u)-(1-\theta)\mathfrak{F}(U_{h}^{n},u_{H}^{n})-\theta f(U_{h}^{n-1}),v_h)=0,
\end{split}
\end{equation}
and
\begin{equation}\label{2.34}
((\sigma-\Sigma_h)^{n-\theta},w_h)+(\nabla (u-U_h)^{n-\theta},\nabla w_h)=0.
\end{equation}
Take $(v_h,w_h)=((\mathfrak{Q}_hu-U_{h})^{n-\theta},(\mathfrak{Q}_h\sigma-\Sigma_{h}))^{n-\theta})$ in system (\ref{2.33})-(\ref{2.34}), use (\ref{2.26}) and sum for the resulting equations to get
\begin{equation}\begin{split}\label{2.35}
&\Big{(}\mathcal{D}_t(\mathfrak{Q}_hu-U_{h})^{n-\theta},(\mathfrak{Q}_hu-U_{h})^{n-\theta}\Big{)}+\gamma\|(\mathfrak{Q}_h\sigma-\Sigma_h)^{n-\theta}\|^2+\|\nabla (\mathfrak{Q}_hu-U_{h})^{n-\theta}\|^2\\=&-\Big{(}\mathcal{D}_t(u-\mathfrak{Q}_hu)^{n-\theta},(\mathfrak{Q}_hu-U_{h})^{n-\theta}\Big{)}
-\gamma((\sigma-\mathfrak{Q}_h\sigma)^{n-\theta},(\mathfrak{Q}_h\sigma-\Sigma_{h}))^{n-\theta})
\\&-(f^{n-\theta}(u)-(1-\theta)\mathfrak{F}(U_{h}^{n},u_{H}^{n})-\theta f(U_{h}^{n-1}),(\mathfrak{Q}_hu-U_{h})^{n-\theta})
\\=&I+II+III.
\end{split}
\end{equation}
Now we estimate the three terms on the right hand side of (\ref{2.35}). Using Cauchy-Schwarz inequality as well as Young inequality, we make use of
the similar derivation as the one in (\ref{2.48}) to get
\begin{equation}\begin{split}\label{2.36}
I+II=&-\Big{(}\mathcal{D}_t(u-\mathfrak{Q}_hu)^{n-\theta},(\mathfrak{Q}_hu-U_{h})^{n-\theta}\Big{)}
-\gamma((\sigma-\mathfrak{Q}_h\sigma)^{n-\theta},(\mathfrak{Q}_h\sigma-\Sigma_{h}))^{n-\theta})
%\\\leq&
%\frac12\|\mathcal{D}_t(u-\mathfrak{Q}_hu)^{n-\theta}\|^2+\frac12\|(\mathfrak{Q}_hu-U_{h})^{n-\theta}\|^2
%+\frac{\gamma}{2}\|(\sigma-\mathfrak{Q}_h\sigma)^{n-\theta}\|^2+\frac{\gamma}{2}\|(\mathfrak{Q}_h\sigma-\Sigma_{h})^{n-\theta}\|^2
%\\\leq&\frac12\Big{\|}\frac{(3-2\theta)[(u-\mathfrak{Q}_hu)^{n}-(u-\mathfrak{Q}_hu)^{n-1}]
%-(1-2\theta)[(u-\mathfrak{Q}_hu)^{n-1}-(u-\mathfrak{Q}_hu)^{n-2}]}{2\Delta t}\Big{\|}^2\\&+\frac12\|(\mathfrak{Q}_hu-U_{h})^{n-\theta}\|^2
%+\frac{\gamma}{2}\|(\sigma-\mathfrak{Q}_h\sigma)^{n-\theta}\|^2+\frac{\gamma}{2}\|(\mathfrak{Q}_h\sigma-\Sigma_{h})^{n-\theta}\|^2
\\\leq&\frac{C}{\Delta t}\int_{t_{n-2}}^{t_n}\|(u-\mathfrak{Q}_hu)_t\|^2dt+\frac12\|(\mathfrak{Q}_hu-U_{h})^{n-\theta}\|^2
\\&+\frac{\gamma}{2}\|(\sigma-\mathfrak{Q}_h\sigma)^{n-\theta}\|^2+\frac{\gamma}{2}\|(\mathfrak{Q}_h\sigma-\Sigma_{h})^{n-\theta}\|^2.
\end{split}
\end{equation}
Use Taylor formula, Cauchy-Schwarz inequality to get
\begin{equation}\begin{split}\label{2.37}
III=&-(f^{n-\theta}(u)-(1-\theta)\mathfrak{F}(U_{h}^{n},u_{H}^{n})-\theta f(U_{h}^{n-1}),(\mathfrak{Q}_hu-U_{h})^{n-\theta})
\\\leq& \|(1-\theta)[f^{n}-\mathfrak{F}(U_{h}^{n},u_{H}^{n})]+\theta[f^{n-1}- f(U_{h}^{n-1})]\|\|(\mathfrak{Q}_hu-U_{h})^{n-\theta}\|
\\=& \|(1-\theta)[f^{n}-f(u_{H}^{n})-f'(u_{H}^{n})(U_{h}^{n}-u_{H}^{n})]+\theta[f^{n-1}- f(U_{h}^{n-1})]\|\|(\mathfrak{Q}_hu-U_{h})^{n-\theta}\|
\\=& \|(1-\theta)[f'(u_{H}^{n})(u^n-U_{h}^{n})+\frac12f''(\phi_1)(u^{n}-u_{H}^{n})^2]
+\theta[f'(\phi_2)(u^{n-1}-U_{h}^{n-1})]\|\|(\mathfrak{Q}_hu-U_{h})^{n-\theta}\|
\\\leq& C(\|u^n-U_{h}^{n}\|^2+\|(u^{n}-u_{H}^{n})^2\|^2
+\|u^{n-1}-U_{h}^{n-1}\|^2+\|(\mathfrak{Q}_hu-U_{h})^{n-\theta}\|^2).
\end{split}
\end{equation}
Substitute (\ref{2.36}) and (\ref{2.37}) into (\ref{2.35}) and use inequality (\ref{4.15}) to get
\begin{equation}\begin{split}\label{2.38}
&\frac{1}{4\Delta t}(\mathbb{H}[(\mathfrak{Q}_hu-U_{h})^{n}]-\mathbb{H}[(\mathfrak{Q}_hu-U_{h})^{n-1}])
+\frac{\gamma}{2}\|(\mathfrak{Q}_h\sigma-\Sigma_h)^{n-\theta}\|^2+\|\nabla (\mathfrak{Q}_hu-U_{h})^{n-\theta}\|^2
\\\leq&\frac{C}{\Delta t}\int_{t_{n-2}}^{t_n}\|(u-\mathfrak{Q}_hu)_t\|^2dt+\frac12\|(\mathfrak{Q}_hu-U_{h})^{n-\theta}\|^2
+\frac{\gamma}{2}\|(\sigma-\mathfrak{Q}_h\sigma)^{n-\theta}\|^2
\\&+C(\|u^n-U_{h}^{n}\|^2+\|(u^{n}-u_{H}^{n})^2\|^2
+\|u^{n-1}-U_{h}^{n-1}\|^2+\|(\mathfrak{Q}_hu-U_{h})^{n-\theta}\|^2).
\end{split}
\end{equation}
Sum (\ref{2.38}) with respect to $n$ from $2$ to $L$ to get
\begin{equation}\begin{split}\label{2.39}
&\mathbb{H}[(\mathfrak{Q}_hu-U_{h})^{L}]
+2\Delta t\sum_{n=2}^L(\gamma\|(\mathfrak{Q}_h\sigma-\Sigma_h)^{n-\theta}\|^2+2\|\nabla (\mathfrak{Q}_hu-U_{h})^{n-\theta}\|^2)
\\\leq&\mathbb{H}[(\mathfrak{Q}_hu-U_{h})^{1}]+C\int_{t_{0}}^{t_L}\|(u-\mathfrak{Q}_hu)_t\|^2dt+2\Delta t\sum_{n=2}^L\|(\mathfrak{Q}_hu-U_{h})^{n-\theta}\|^2
+2\gamma\Delta t\sum_{n=2}^L\|(\sigma-\mathfrak{Q}_h\sigma)^{n-\theta}\|^2
\\&+C\Delta t\sum_{n=2}^L(\|u^n-U_{h}^{n}\|^2+\|(u^{n}-u_{H}^{n})^2\|^2
+\|u^{n-1}-U_{h}^{n-1}\|^2+\|(\mathfrak{Q}_hu-U_{h})^{n-\theta}\|^2).
\end{split}
\end{equation}
Now we need to give the estimate of $\mathbb{H}[(\mathfrak{Q}_hu-U_{h})^{1}]$. We combine (\ref{2.21})-(\ref{2.31}) with (\ref{2.151})-(\ref{2.161}) to get
\begin{equation}\begin{split}\label{2.405}
\Big{(}\mathcal{D}_t(u-U_h)^{\frac12},v_h\Big{)}&-\gamma(\nabla (\sigma-\Sigma_h)^{\frac12},\nabla v_h)+(\nabla (u-U_h)^{\frac12},\nabla v_h)\\&+(f^{\frac12}(u)-\frac12\mathfrak{F}(U_{h}^{1},u_{H}^{1})-\frac12 f(U_{h}^{0}),v_h)=0,
\end{split}
\end{equation}
and
\begin{equation}\label{2.415}
((\sigma-\Sigma_h)^{\frac12},w_h)+(\nabla (u-U_h)^{\frac12},\nabla w_h)=0.
\end{equation}
In (\ref{2.405})-(\ref{2.415}), we take $(v_h,w_h)=((\mathfrak{Q}_hu-U_{h})^{\frac12},(\mathfrak{Q}_h\sigma-\Sigma_{h})^{\frac12})$ and use the similar derivation to the one of inequality (\ref{2.39}) to get
\begin{equation}\begin{split}\label{2.42}
&\|(\mathfrak{Q}_hu-U_{h})^{1}\|^2
+2\Delta t(\gamma\|(\mathfrak{Q}_h\sigma-\Sigma_h)^{\frac12}\|^2+2\|\nabla (\mathfrak{Q}_hu-U_{h})^{\frac12}\|^2)
\\\leq&\|(\mathfrak{Q}_hu-U_{h})^{0}\|^2+C\int_{t_{0}}^{t_1}\|(u-\mathfrak{Q}_hu)_t\|^2dt+2\Delta t\|(\mathfrak{Q}_hu-U_{h})^{\frac12}\|^2
+2\gamma\Delta t\|(\sigma-\mathfrak{Q}_h\sigma)^{\frac12}\|^2
\\&+C\Delta t(\|u^1-U_{h}^{1}\|^2+\|(u^{1}-u_{H}^{1})^2\|^2
+\|u^{0}-U_{h}^{0}\|^2+\|(\mathfrak{Q}_hu-U_{h})^{\frac12}\|^2).
\end{split}
\end{equation}
By the similar process to the derivation (\ref{2.55}), we have
\begin{equation}\begin{split}\label{2.43}
\mathbb{H}[(\mathfrak{Q}_hu-U_{h})^{1}]\leq & C(\|(\mathfrak{Q}_hu-U_{h})^{1}\|^2+\|(\mathfrak{Q}_hu-U_{h})^{0}\|^2).
\end{split}
\end{equation}
Combining (\ref{2.42}), (\ref{2.43}) with (\ref{2.39}) and using Gronwall lemma with (\ref{5.17}), we have
\begin{equation}\begin{split}\label{2.44}
&\|(\mathfrak{Q}_hu-U_{h})^{L}\|^2
+2\Delta t\sum_{n=1}^L(\gamma\|(\mathfrak{Q}_h\sigma-\Sigma_h)^{n-\theta}\|^2+2\|\nabla (\mathfrak{Q}_hu-U_{h})^{n-\theta}\|^2)
\\\leq&C\int_{t_{0}}^{t_L}\|(u-\mathfrak{Q}_hu)_t\|^2dt+C\Delta t\sum_{n=1}^L(\|(u^{n}-u_{H}^{n})^2\|^2
+h^{2m+2})+C\Delta th^{2m+2}
\\\leq&C\int_{t_{0}}^{t_L}\|(u-\mathfrak{Q}_hu)_t\|^2dt+C\Delta t\sum_{n=1}^L(H^{4m+4}
+h^{2m+2})+C\Delta th^{2m+2}.
\end{split}
\end{equation}
Using triangle inequality, we arrive at the conclusion (\ref{5.18}).
\par
Combine (\ref{1.1})-(\ref{1.2}), (\ref{1.21})-(\ref{1.22}) with the above results to arrive at the main error theorem in this paper.
\begin{thm}\label{thm6}There exists a constant $C$ independent of space-time mesh step length $h,H,\Delta t$ such that
\begin{equation}\label{2.6}
\|u(t_n)-u_h^{n}\|+\gamma^{\frac12}\Big{(}\Delta t\sum_{k=1}^n\|\sigma(t_{k-\theta})-\sigma_h^{k-\theta}\|^2\Big{)}^{\frac12}+\Big{(}\Delta t\sum_{k=1}^n\|\nabla (u(t_{k-\theta})-u_{h}^{k-\theta})\|^2\Big{)}^{\frac12}\leq C[\Delta t^2+h^{m+1}],
\end{equation}
\begin{equation}\label{2.61}
\|u(t_n)-u_H^{n}\|+\gamma^{\frac12}\Big{(}\Delta t\sum_{k=1}^n\|\sigma(t_{k-\theta})-\sigma_H^{k-\theta}\|^2\Big{)}^{\frac12}+\Big{(}\Delta t\sum_{k=1}^n\|\nabla (u(t_{k-\theta})-u_{H}^{k-\theta})\|^2\Big{)}^{\frac12}\leq C[\Delta t^2+H^{m+1}],
\end{equation}
and
\begin{equation}\label{2.62}
\|u(t_n)-U_h^{n}\|+\gamma^{\frac12}\Big{(}\Delta t\sum_{k=1}^n\|\sigma(t_{k-\theta})-\Sigma_h^{k-\theta}\|^2\Big{)}^{\frac12}+\Big{(}\Delta t\sum_{k=1}^n\|\nabla (u(t_{k-\theta})-U_{h}^{k-\theta})\|^2\Big{)}^{\frac12}\leq C[\Delta t^2+h^{m+1}+H^{2m+2}].
\end{equation}
\end{thm}
\textbf{Proof.} Let $E_u^{n-\theta}=u(t_{n-\theta})-u^{n-\theta}$ and $E_\sigma^{n-\theta}=\sigma(t_{n-\theta})-\sigma^{n-\theta}$
and we easily know from systems (\ref{1.1})-(\ref{1.2}) and (\ref{1.12})-(\ref{1.23}) that for $n\geq2$
\begin{equation}\begin{split}\label{5.2}
\Big{(}\mathcal{D}_tE_u^{n-\theta},v\Big{)}-\gamma(\nabla E_\sigma^{n-\theta},\nabla v)+(\nabla E_u^{n-\theta},\nabla v)+(f(u(t_{n-\theta}))-f^{n-\theta}(u),v)=&(R_1^{n-\theta},v),~\forall v\in H_0^1,
\end{split}
\end{equation}
and
\begin{equation}\label{5.3}
(E_\sigma^{n-\theta},w)+(\nabla E_u^{n-\theta},\nabla w)=(R_2^{n-\theta},w),~\forall w\in H_0^1,
\end{equation}
where
\begin{equation}\begin{split}\label{5.4}R_1^{n-\theta}=&\mathcal{D}_t(u(t_{n-\theta})-u^{n-\theta})+\gamma\bigtriangleup (\sigma(t_{n-\theta})-\sigma^{n-\theta})-\bigtriangleup (u(t_{n-\theta})-u^{n-\theta})+(f(u(t_{n-\theta}))-f^{n-\theta}(u))\\=&
O(\Delta t^2),\end{split}
\end{equation}
\begin{equation}\begin{split}\label{5.5}R_2^{n-\theta}=&(\sigma(t_{n-\theta})-\sigma^{n-\theta}))-\bigtriangleup (u(t_{n-\theta})-u^{n-\theta})=
O(\Delta t^2),\end{split}
\end{equation}
Also, we easily follow that
\begin{equation}\begin{split}\label{5.6}
\Big{(}\mathcal{D}_tE_u^{\frac12},v\Big{)}-\gamma(\nabla E_\sigma^{\frac12},\nabla v)+(\nabla E_u^{\frac12},\nabla v)+(f(u(t_{\frac12}))-f^{\frac12}(u),v)=&(R_3^{\frac12},v),~\forall v\in H_0^1,
\end{split}
\end{equation}
and
\begin{equation}\label{5.7}
(E_\sigma^{\frac12},w)+(\nabla E_u^{\frac12},\nabla w)=(R_4^{\frac12},w),~\forall w\in H_0^1,
\end{equation}
where
\begin{equation}\begin{split}\label{5.8}R_3^{\frac12}=&\mathcal{D}_t(u(t_{\frac12})-u^{\frac12})+\gamma\bigtriangleup (\sigma(t_{\frac12})-\sigma^{\frac12})-\bigtriangleup (u(t_{n-\theta})-u^{n-\theta})+(f(u(t_{\frac12}))-f^{\frac12}(u))\\=&
O(\Delta t^2),\end{split}
\end{equation}
\begin{equation}\begin{split}\label{5.9}R_4^{\frac12}=&(\sigma(t_{\frac12})-\sigma^{\frac12}))-\bigtriangleup (u(t_{\frac12})-u^{\frac12})=
O(\Delta t^2).\end{split}
\end{equation}
We use the similar analysis as the ones in the theorem \ref{thm5} to easily get
 \begin{equation}\label{5.10}
\|u(t_n)-u^n\|+\gamma^{\frac12}\Big{(}\Delta t\sum_{k=1}^n\|\sigma(t_{k-\theta})-\sigma^{k-\theta}\|^2\Big{)}^{\frac12}+\Big{(}\Delta t\sum_{k=1}^n\|\nabla u(t_{k-\theta})-u^{k-\theta}\|^2\Big{)}^{\frac12}\leq C\Delta t^2.
\end{equation}
We combine (\ref{5.10}) with the results in theorem \ref{thm5} and use triangle inequality to get the conclusions of theorem \ref{thm6}.
\section{Numerical tests}
In this section, we take some numerical examples to test the computational efficiency of TGMFE method combined with $\theta$-scheme with second-order convergence rate. For implementing the numerical computations in two-dimensional cases, we take rectangular partition for spatial domain $\overline{\Omega}$ and choose continuous bilinear element with basis function $P(x_1,x_2)=a+bx_1+cx_2+dx_1x_2$. We also consider a one-dimensional case based on the piecewise continuous linear element. In the following discussions, we take $\hat{H}$ as the length of the edge of coarse grid rectangular unit and $\hat{h}$ the length of the edge of fine grid rectangular unit, respectively. It is easy to see that we have $H=\sqrt{2}\hat{H}$ and $h=\sqrt{2}\hat{h}$.
%\begin{figure}[h]
%\begin{center}
%\begin{minipage}{15cm}
%  % Requires \usepackage{graphicx}
%  \centering\includegraphics[width=15cm]{rec_grid.eps}
%  \caption{The type of meshes}\label{C0}
%\end{minipage}
%\end{center}
%\end{figure}
\par
\textbf{Example 4.1}
\\
In (\ref{0.1}), based on the spatial domain $\bar{\Omega}=[-1,1]\times[-1,1]$ and the temporal interval $[0,1]$, we take nonlinear term $f(u)=u^3-u$ and the exact solution
$u(\textbf{z},t)=\exp(-t)\sin(2\pi x_1)\sin(2\pi x_2),\textbf{z}=(x_1,x_2)$,
then we get the source term $g(\textbf{z},t)=(8\pi^2-2+64\gamma\pi^4)\exp(-t)\sin(2\pi x_1)\sin(2\pi x_2)
 +\exp(-3t)\sin^3(2\pi x_1)\sin^3(2\pi x_2)$.
 \par
 In Table \ref{tab1}, by taking $\theta=0.2$, changed $\gamma=0.01,1,10$ and $\Delta t=\hat{h}=\hat{H}^2=1/25,1/64,1/100$, we arrive at TGMFE and nonlinear Galerkin MFE error estimates with second-order convergence rate. By the comparison between TGMFE method and standard nonlinear Galerkin MFE method, ones see that TGMFE method can save the CPU time largely. We also see from Tables \ref{tab2}-\ref{tab4} that the similar results can be obtained based on the cases $\theta=0.4,0,0.5$. Compared with nonlinear MFE method in this example, ones can see clearly from the calculated data in Tables \ref{tab2}-\ref{tab4} that the TGMFE method not only save the computing time, but also get the better convergence rate.
 \begin{table}[h]
\centering
 \caption{$\theta=0.2$, $H=\sqrt{2}\hat{H}$ and  $h=\sqrt{2}\hat{h}$}\label{tab1}
 \begin{tabular}{*{10}{c}}
 \midrule
   %\bottomrule
 %  &&&&&&\\
$\gamma$&$\Delta t$& $\hat{H}$&$\hat{h}$&	$\|u-U_h\|$	&Order&$\|\sigma-\Sigma_h\|$& Order&	$\frac{\text{CPU time}}{\text{Seconds}}$ \\
\midrule
0.01&1/25& 	   1/5 & 	   1/25 &	4.41472E-02&	&	2.06945E+00&	&	0.3545\\
   &1/64 &	   1/8 & 	   1/64 &	5.03765E-03&	2.30912 &	1.97903E-01&	2.49706 	&5.3805\\
   &1/100&	   1/10& 	   1/100&	1.88823E-03&	2.19881 	&7.27048E-02&	2.24378 	&21.6419\\						
  1& 1/25& 	   1/5&  	   1/25& 	8.05926E-02&	&	4.91240E+00&	&	0.3970\\
  & 1/64& 	   1/8&  	   1/64& 	1.25046E-02&	1.98223& 	7.44673E-01&	2.00698 &	5.3017\\
  & 1/100&	   1/10 &	   1/100&	5.08477E-03&	2.01629 &	3.00895E-01&	2.03050& 	21.5775\\
		
  10& 1/25 &	   1/5 & 	   1/25& 	8.16219E-02&&		4.99464E+00&&		0.3514\\
  & 1/64& 	   1/8 & 	   1/64 &	1.28273E-02&	1.96862 &	7.69990E-01&	1.98907 &	5.3197\\
  & 1/100&	   1/10& 	   1/100&	5.27285E-03&	1.99200 &	3.15616E-01&	1.99838 &	21.4013\\
\hline
MFE method&&&&$\|u-u_h\|$&&$\|\sigma-\sigma_h\|$&&	\\
\hline							
  0.01 &1/25& 		&   1/25 &	4.58575E-02&	&	2.19476E+00&	&	0.8032\\
   &1/64& 		&   1/64& 	5.45864E-03&	2.26417& 	2.21379E-01&	2.44035& 	15.1652\\
   &1/100&		&   1/100&	2.10826E-03&	2.13167& 	8.27276E-02&	2.20559 &	63.9630\\
									
1  & 1/25 	&	&   1/25 &	8.06655E-02	&&	4.91823E+00&&		0.7560\\
  & 1/64 	&	&   1/64 &	1.25411E-02&	1.98009& 	7.47534E-01&	2.00416 &	15.0554\\
  & 1/100	&	&   1/100&	5.10577E-03&	2.01359& 	3.02537E-01&	2.02689 &	63.2929\\
											
 10 & 1/25 &		&   1/25& 	8.16294E-02&&		4.99523E+00&&		0.7677\\
  & 1/64 &		&   1/64 &	1.28312E-02&	1.96840 &	7.70294E-01&	1.98878& 	15.3985\\
  & 1/100&		&   1/100&	5.27514E-03&	1.99171 &	3.15795E-01&	1.99799 &	65.6151\\
 \bottomrule
 \end{tabular}
 \end{table}

\begin{table}[h]
\centering
 \caption{$\theta=0.4$,$H=\sqrt{2}\hat{H}$ and  $h=\sqrt{2}\hat{h}$}\label{tab2}
 \begin{tabular}{*{10}{c}}
 \midrule
   %\bottomrule
 %  &&&&&&\\
$\gamma$&$\Delta t$& $\hat{H}$&$\hat{h}$&	$\|u-U_h\|$	&Order&$\|\sigma-\Sigma_h\|$& Order&	$\frac{\text{CPU time}}{\text{Seconds}}$\\
\midrule
  1   &1/25& 	  1/5  & 	   1/25	&	8.05926E-02	&&	4.91240E+00	&&	0.3353\\
   &1/64 &	   1/8 & 	   1/64 &	1.25046E-02&	1.98223 &	7.44673E-01	&2.00698& 	5.2955\\
  &1/100&	   1/10& 	   1/100&	5.08477E-03&	2.01629& 	3.00895E-01&	2.03050& 	21.4680\\
\hline
MFE method&&&&$\|u-u_h\|$&&$\|\sigma-\sigma_h\|$&&	\\
\hline									
   &1/25& 	    & 	   1/25	&	8.06655E-02	&&	4.91823E+00	&&	0.7715\\
  &1/64 &	    & 	   1/64 &	1.25411E-02&	1.98009 &	7.47534E-01&	2.00416 &	15.2819\\
   &1/100&	   & 	   1/100&	5.10577E-03&	2.01359 	&3.02537E-01	&2.02689 &	65.1652\\
 \bottomrule
 \end{tabular}
 \end{table}

\begin{table}[h]
\centering
 \caption{$\theta=0$, $H=\sqrt{2}\hat{H}$ and  $h=\sqrt{2}\hat{h}$}\label{tab3}
 \begin{tabular}{*{10}{c}}
 \midrule
   %\bottomrule
 %  &&&&&&\\
$\gamma$&$\Delta t$& $\hat{H}$&$\hat{h}$&	$\|u-U_h\|$	&Order&$\|\sigma-\Sigma_h\|$& Order&	$\frac{\text{CPU time}}{\text{Seconds}}$\\
\midrule
	   1&1/25& 	   1/5 & 	   1/25&	8.05926E-02	&&	4.91240E+00&&		0.3515\\
    &1/64 &	   1/8 & 	   1/64 & 	1.25046E-02&	1.98223 &	7.44673E-01	&2.00698 &	5.3909\\
  &1/100&	   1/10& 	   1/100&	5.08477E-03&	2.01629& 	3.00895E-01&	2.03050& 	21.2432\\
\hline
MFE method&&&&$\|u-u_h\|$&&$\|\sigma-\sigma_h\|$&&	\\
\hline									
   &1/25& 	    & 	   1/25	&8.06655E-02		&&4.91823E+00		&&0.7599\\
    &1/64 &	    & 	   1/64 &	1.25411E-02&	1.98009 &	7.47534E-01&	2.00416 &	15.2091\\
   &1/100&	   & 	   1/100&	5.10577E-03&	2.01359 &	3.02537E-01&	2.02689 	&65.2895\\						
 \bottomrule
 \end{tabular}
 \end{table}

\begin{table}[h]
\centering
 \caption{$\theta=0.5$, $H=\sqrt{2}\hat{H}$ and  $h=\sqrt{2}\hat{h}$}\label{tab4}
 \begin{tabular}{*{10}{c}}
 \midrule
   %\bottomrule
 %  &&&&&&\\
$\gamma$&$\Delta t$& $\hat{H}$&$\hat{h}$&	$\|u-U_h\|$	&Order&$\|\sigma-\Sigma_h\|$& Order&	$\frac{\text{CPU time}}{\text{Seconds}}$\\
\midrule
      1   &1/25& 	  1/5  & 	   1/25	&	8.05926E-02	&&	4.91240E+00&	&	0.3523\\
   &1/64 &	   1/8 & 	   1/64 &	1.25046E-02	&1.98223 &	7.44673E-01	&2.00698 &	5.3813\\
   &1/100&	   1/10& 	   1/100&	5.08477E-03&	2.01629 &	3.00895E-01&	2.03050& 	21.7946\\
\hline
MFE method&&&&$\|u-u_h\|$&&$\|\sigma-\sigma_h\|$&&	\\
\hline									
   &1/25& 	    & 	   1/25	& 	8.06655E-02	&&	4.91823E+00	&&	0.7726\\
   &1/64 &	    & 	   1/64 &	1.25411E-02	&1.98009 &	7.47534E-01&	2.00416& 	15.1932\\
  &1/100&	   & 	   1/100&	5.10577E-03&	2.01359 	&3.02537E-01&	2.02689 &	63.2262\\
 \bottomrule
 \end{tabular}
 \end{table}
\textbf{Example 4.2}
\\
Based on the spatial domain $\bar{\Omega}=[0,1]\times[0,1]$ and the time interval $[0,1]$ and taking the source term $g(\textbf{z},t)=0$, the nonlinear term $f(u)=u^3-u$, we choose another numerical example only including the initial value $u(\textbf{z},0)=x^3(1-x)^3y^3(1-y)^3$ to verify the convergence order in space.
\par
For testing the spatial convergence order, we choose the numerical solution pair $\{U_h,\Sigma_h\}$ with $\hat{h}=1/100$ and $\Delta t=1/200$ as the approximate exact solution. In Tables \ref{tab5}-\ref{tab8}, with the parameter $\gamma=0.1$, the fixed time step length $\Delta t=1/200$ and changed $\theta=0,0.1,0.3,0.5$, we make the errors between the case $\hat{h}=\hat{H}^2=1/100$ with the cases $\hat{h}=\hat{H}^2=1/9,1/16,1/25$ to get the convergence results, which show the convergence order of our method is close to $2$. Similarly, we also arrive at second-order convergence rate based on the case $\gamma=20$ in Tables \ref{tab5}-\ref{tab8}. In the second example, by using TGMFE method, we get almost the same convergence results as the ones of nonlinear MFE method.
\par
\begin{figure}[h]
\begin{center}
\begin{minipage}{7.5cm}
  \centering\includegraphics[width=7cm]{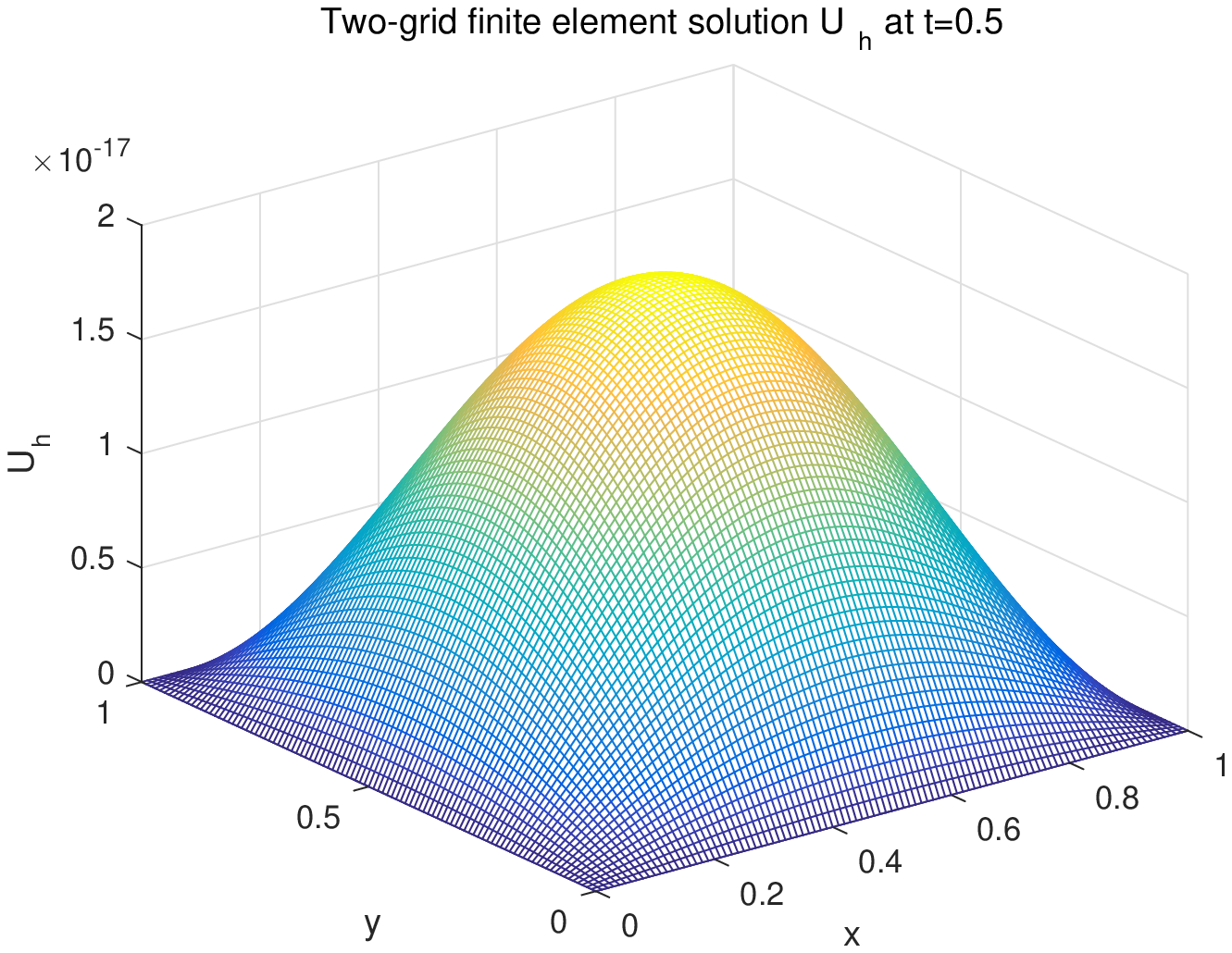}
  \caption{$U_h$ with $\hat{h}=\hat{H}^2=\frac{1}{100}$, $\Delta t=\frac{1}{200}$, $\theta=0.1$ and $\gamma=0.1$}\label{C1}
\end{minipage}
\begin{minipage}{7.5cm}
  \centering\includegraphics[width=7cm]{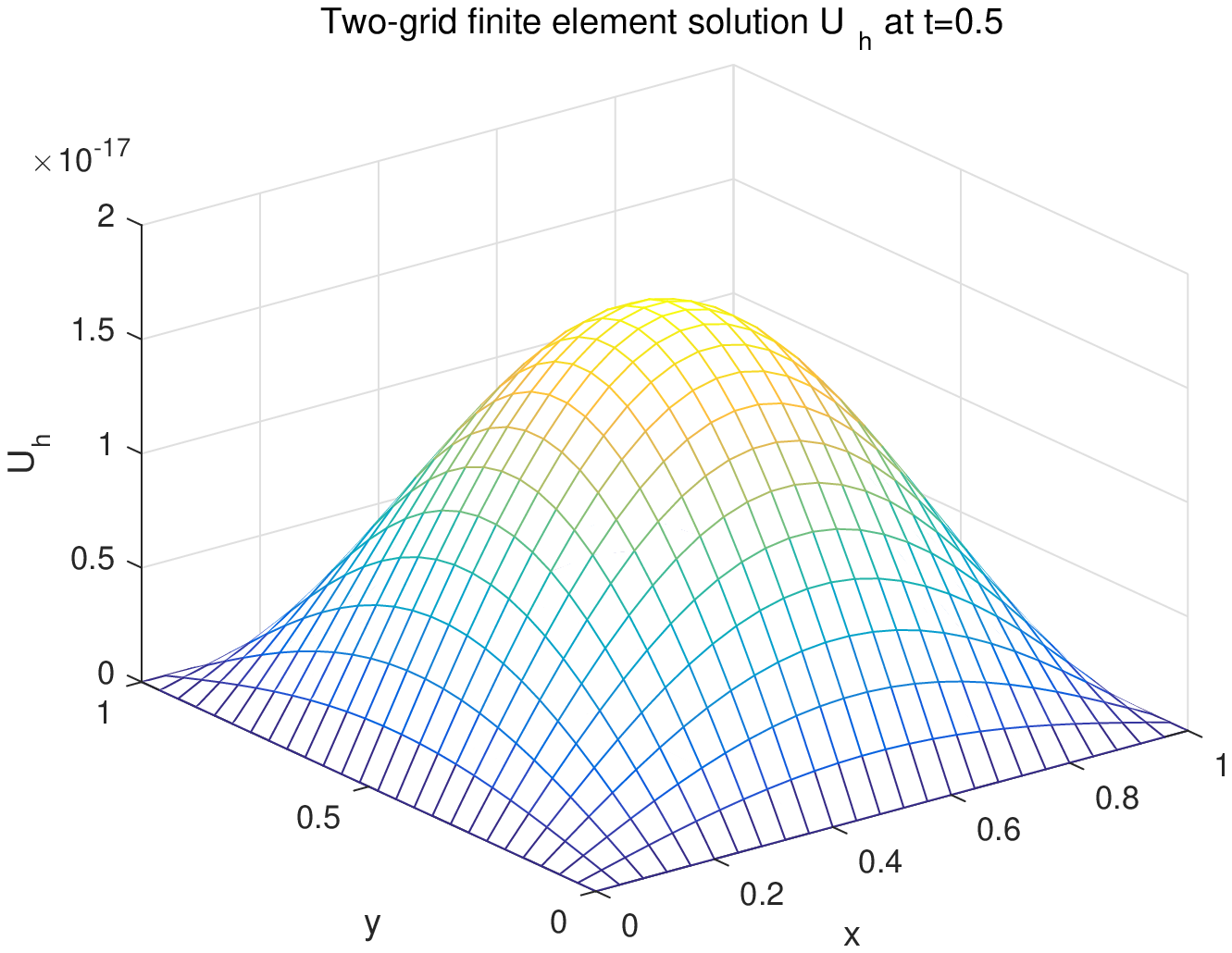}
  \caption{$U_h$ with $\hat{h}=\hat{H}^2=\frac{1}{25}$, $\Delta t=\frac{1}{200}$, $\theta=0.1$ and $\gamma=0.1$}\label{C2}
\end{minipage}%\\[1em]
\end{center}
\end{figure}
\begin{figure}[h]
\begin{center}
\begin{minipage}{7.5cm}
  \centering\includegraphics[width=7cm]{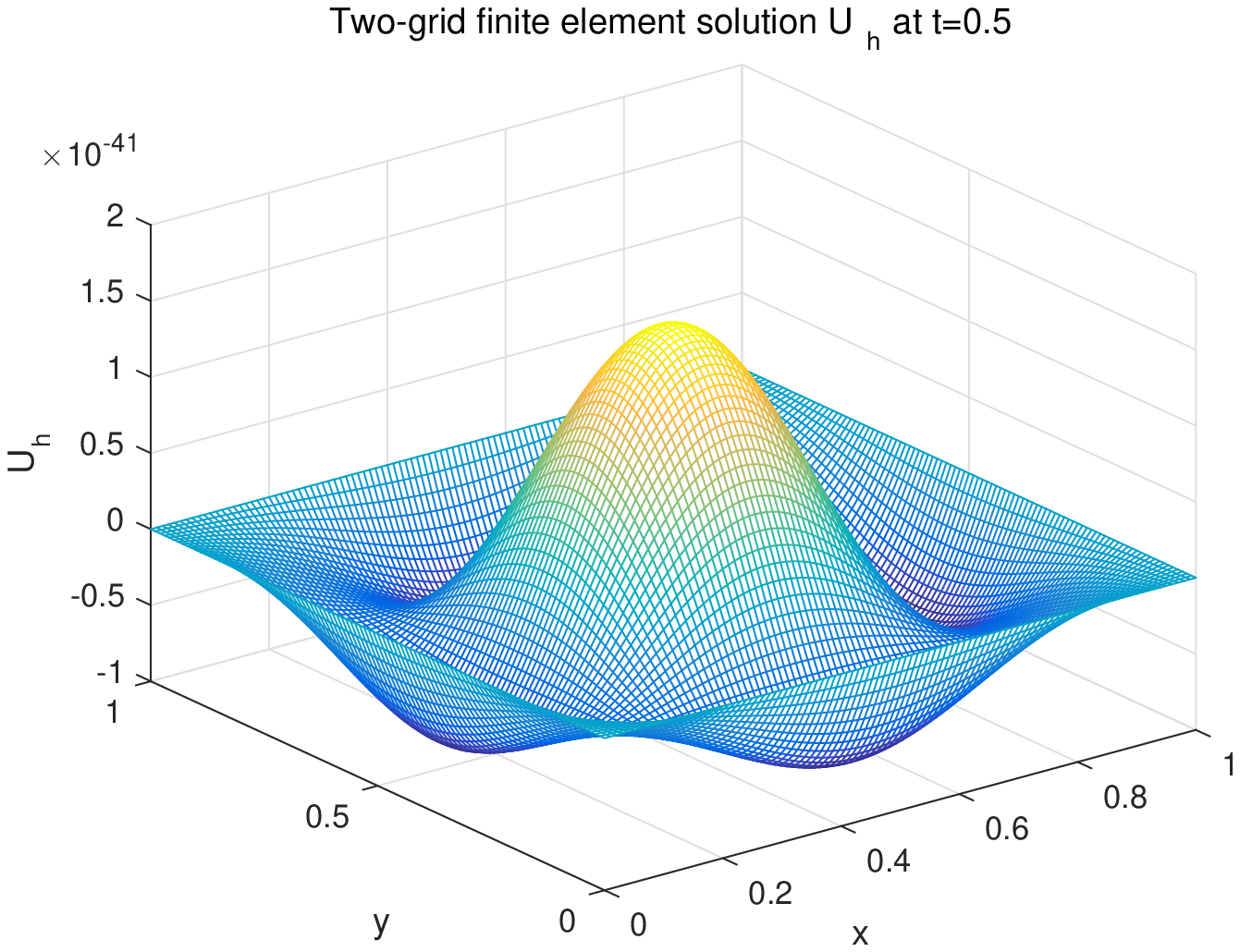}
  \caption{$U_h$ with $\hat{h}=\hat{H}^2=\frac{1}{100}$, $\Delta t=\frac{1}{200}$, $\theta=0.3$ and $\gamma=20$}\label{C3}
\end{minipage}
\begin{minipage}{7.5cm}
  \centering\includegraphics[width=7cm]{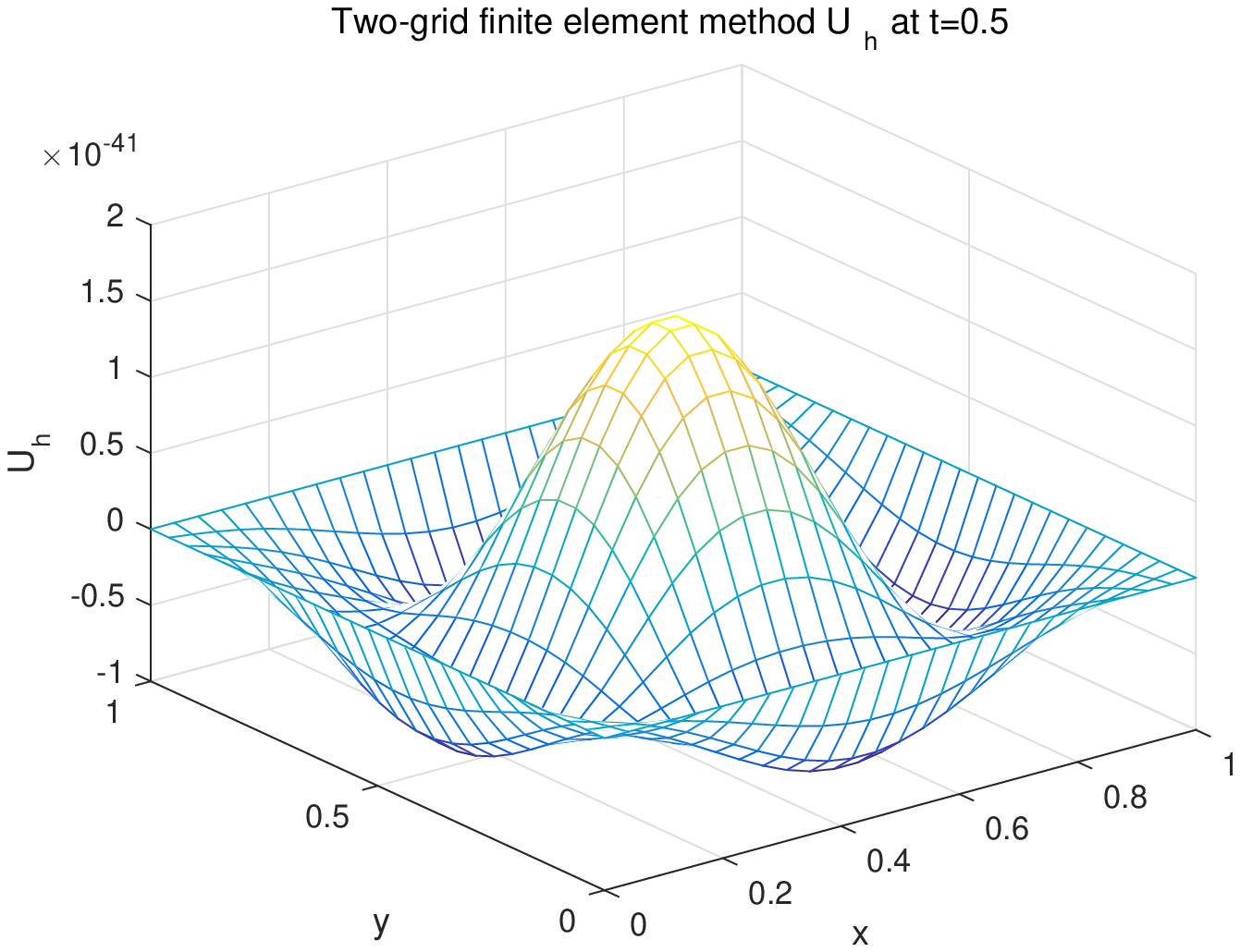}
  \caption{$U_h$ with $\hat{h}=\hat{H}^2=\frac{1}{25}$, $\Delta t=\frac{1}{200}$, $\theta=0.3$ and $\gamma=20$}\label{C4}
\end{minipage}%\\[1em]
\end{center}
\end{figure}
\begin{figure}[h]
\begin{center}
\begin{minipage}{7.5cm}
  \centering\includegraphics[width=7cm]{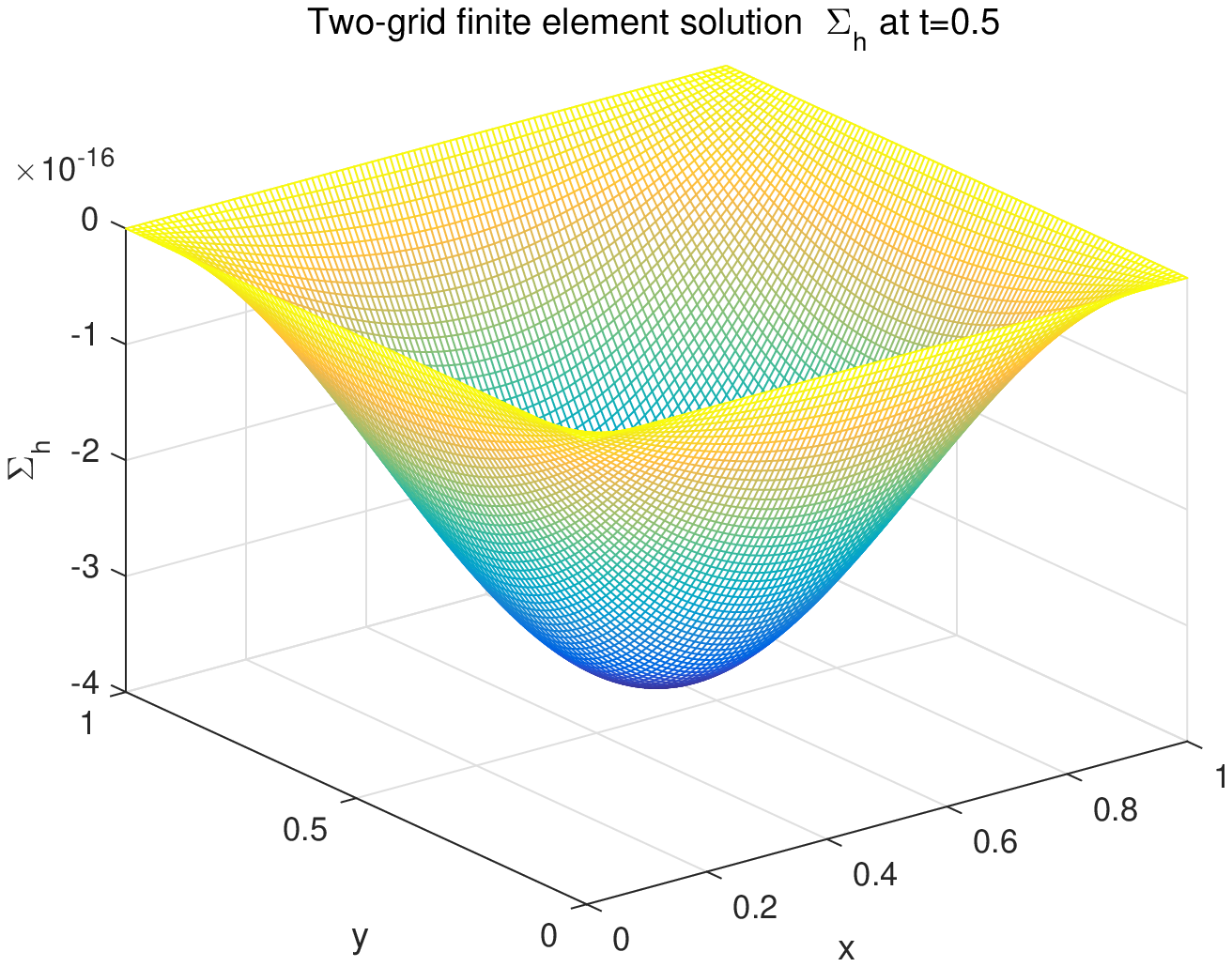}
  \caption{$\Sigma_h$ with $\hat{h}=\hat{H}^2=\frac{1}{100}$, $\Delta t=\frac{1}{200}$, $\theta=0.1$ and $\gamma=0.1$}\label{C5}
\end{minipage}
\begin{minipage}{7.5cm}
  \centering\includegraphics[width=7cm]{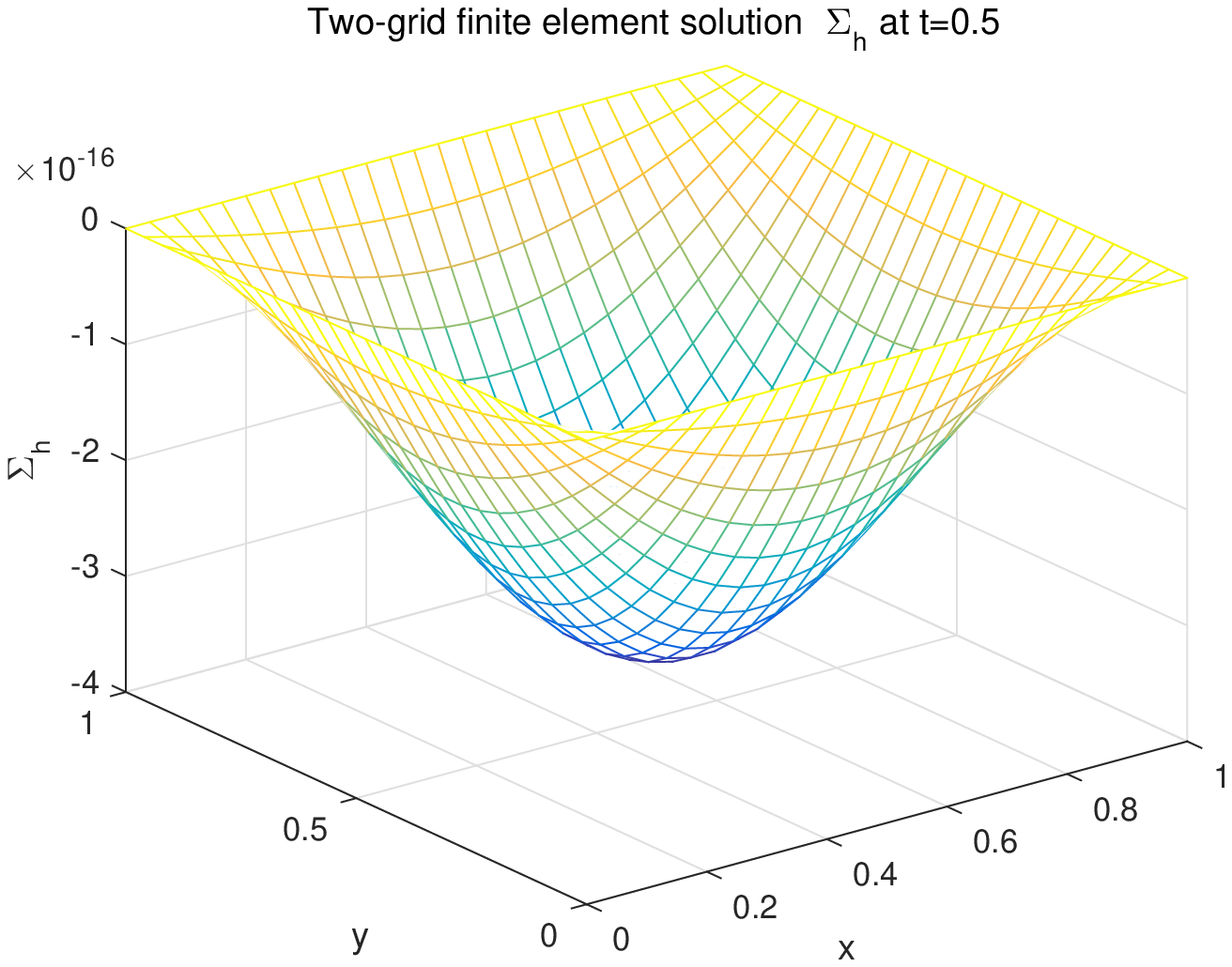}
  \caption{$\Sigma_h$ with $\hat{h}=\hat{H}^2=\frac{1}{25}$, $\Delta t=\frac{1}{200}$, $\theta=0.1$ and $\gamma=0.1$}\label{C6}
\end{minipage}%\\[1em]
\end{center}
\end{figure}
\begin{figure}[h]
\begin{center}
\begin{minipage}{7.5cm}
  \centering\includegraphics[width=7cm]{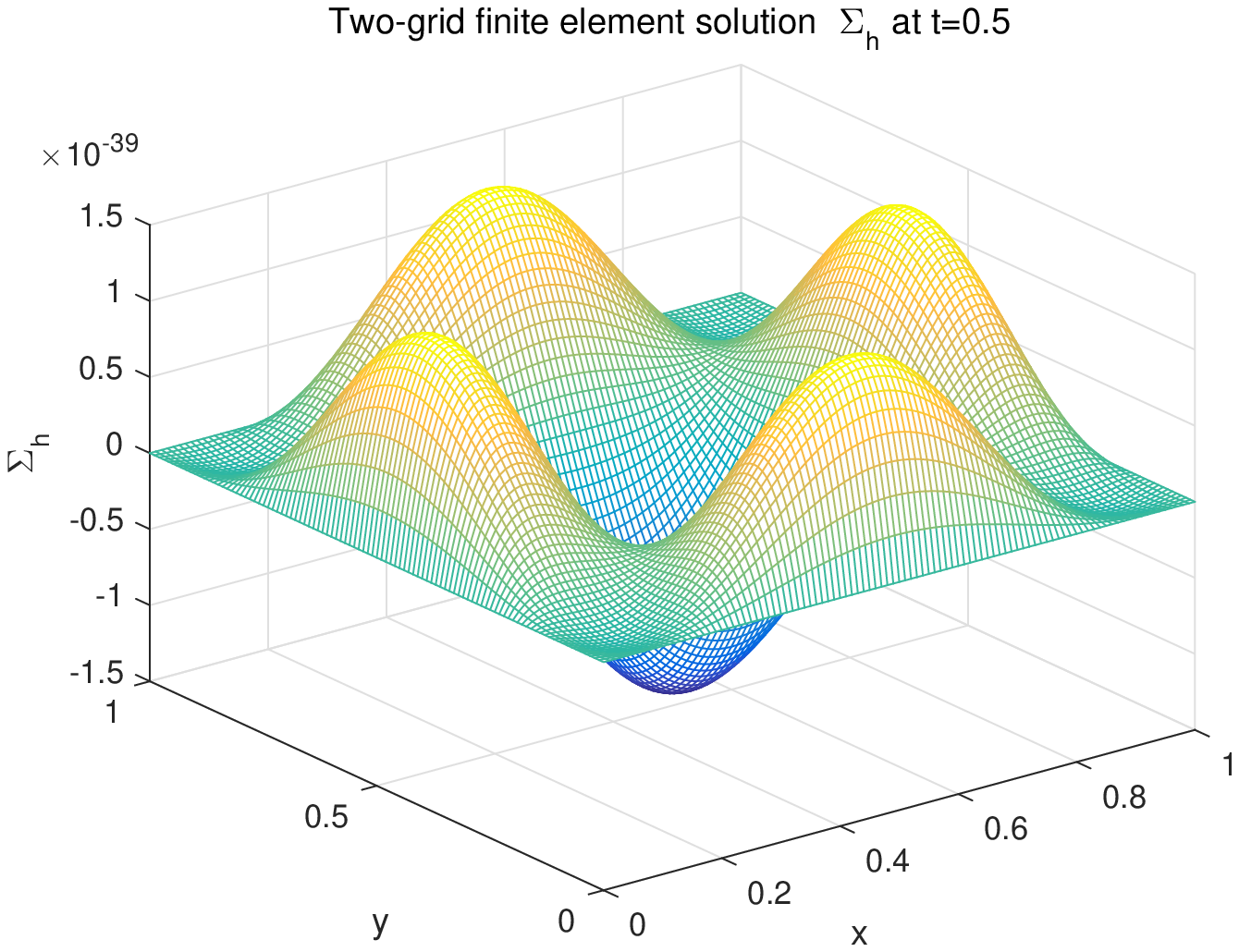}
  \caption{$\Sigma_h$ with $\hat{h}=\hat{H}^2=\frac{1}{100}$, $\Delta t=\frac{1}{200}$, $\theta=0.3$ and $\gamma=20$}\label{C7}
\end{minipage}
\begin{minipage}{7.5cm}
  \centering\includegraphics[width=7cm]{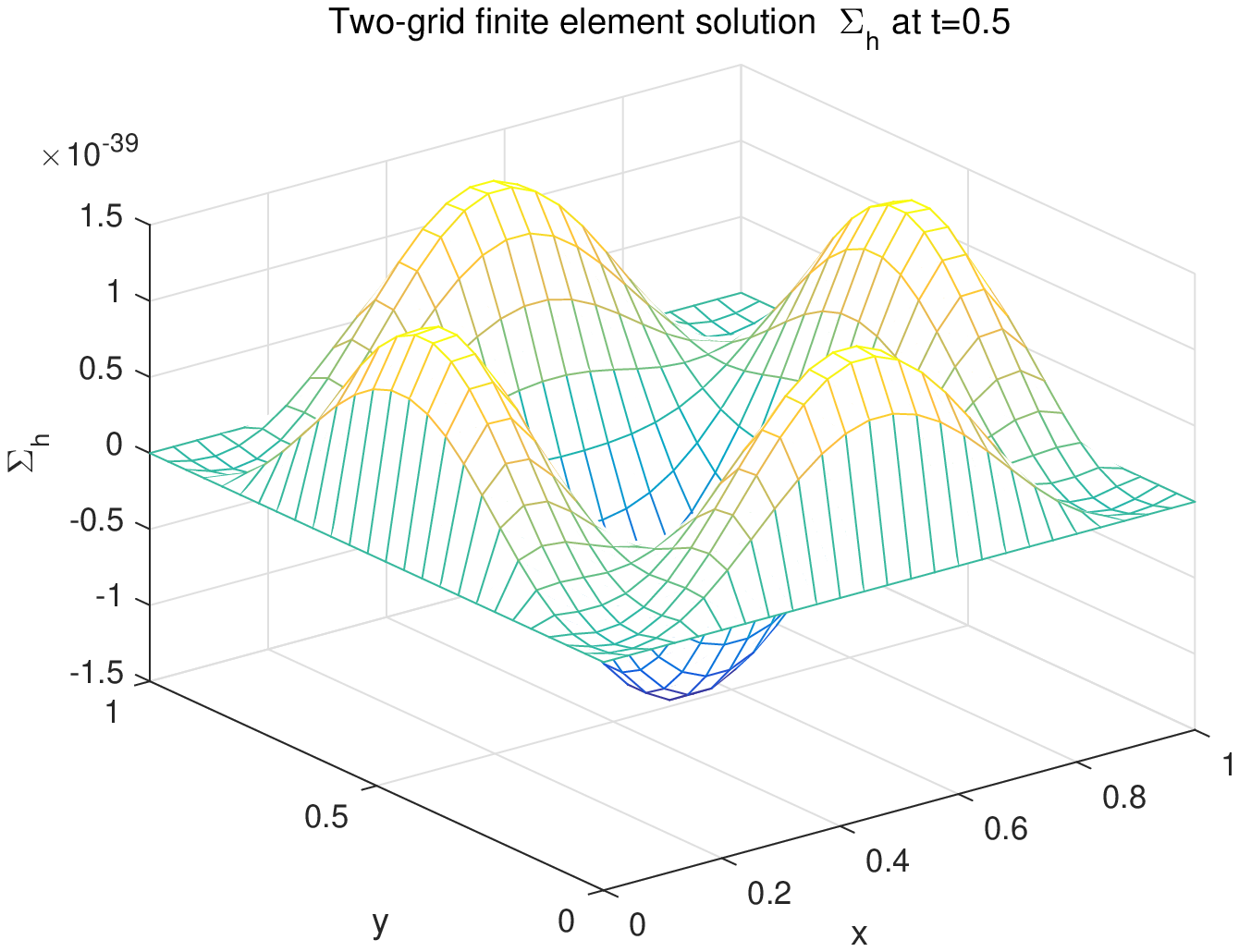}
  \caption{$\Sigma_h$ with $\hat{h}=\hat{H}^2=\frac{1}{25}$, $\Delta t=\frac{1}{200}$, $\theta=0.3$ and $\gamma=20$}\label{C8}
\end{minipage}%\\[1em]
\end{center}
\end{figure}

Further, we show the figures of numerical solutions. In Figs. \ref{C1}-\ref{C2}, based on the parameters $\theta=0.1$, $\gamma=0.1$ and $\Delta t=1/200$, we give the
figures of numerical solution $U_h$ at $t=0.5$ under the parameters $\hat{h}=\hat{H}^2=1/100$ and $\hat{h}=\hat{H}^2=1/25$, respectively. Based on this case, we see that the numerical solutions have the same approximate trend.
In Figs. \ref{C3}-\ref{C4}, choosing the same parameters as in Figs. \ref{C1}-\ref{C2}, we also draw the numerical solutions $U_h$ with $\gamma=20$ and $\theta=0.3$ at $t=0.5$. From Figs. \ref{C3}-\ref{C4}, we can find that the two figures also have the same numerical behavior. By the comparisons in \ref{C1}-\ref{C4}, we see that the numerical solutions $U_h$ with parameters $\gamma=0.1$ and $\gamma=10$ have different numerical behavior. At the same time, we also provide the behavior of numerical solution $\Sigma_h$ in Figs. \ref{C5}-\ref{C8}.
\begin{table}[h]
\centering
 \caption{$\theta=0$, $\Delta t=1/200$, $H=\sqrt{2}\hat{H}$ and  $h=\sqrt{2}\hat{h}$}\label{tab5}
 \begin{tabular}{*{9}{c}}
 \midrule
   %\bottomrule
 %  &&&&&&\\
$\gamma$& $\hat{H}$&$\hat{h}$&	$\|u-U_h\|$	&Order&$\|\sigma-\Sigma_h\|$& Order&	$\frac{\text{CPU time}}{\text{Seconds}}$\\
\midrule
$0.1$      &1/3  	  & 1/9     &2.06411E-06	&	        &4.70651E-04	&	        &0.350528\\
	       &1/4  	  & 1/16 	&6.46351E-07	&2.01805 	&1.49109E-04	&1.99777 	&0.799282\\
	       &1/5  	  & 1/25 	&2.56619E-07	&2.06986 	&5.89125E-05	&2.08077 	&2.055724\\
$20$       &1/3  	  & 1/9  	&1.64081E-06	&	        &4.55079E-04	&	        &0.351065\\
           &1/4  	  & 1/16 	&4.68907E-07	&2.17696 	&1.42591E-04	&2.01696 	&0.792415\\
           &1/5  	  & 1/25 	&1.98760E-07	&1.92321 	&5.64742E-05	&2.07534 	&2.073567\\
   \hline
MFE method&&&$\|u-u_h\|$&&$\|\sigma-\sigma_h\|$&&&	\\
\hline		

		0.1               &&1/9  	&2.06411E-06	&	        &4.70651E-04	&	        &0.484523\\
		              & &1/16 	&6.46351E-07	&2.01805 	&1.49109E-04	&1.99777    &1.674808\\
		              & &1/25   &2.56619E-07	&2.06986    &5.89125E-05	&2.08077 	&5.390575\\

	    20                 & &1/9  	&1.64081E-06	&	        &4.55079E-04	&	        &0.471047\\
	                  & &1/16 	&4.68907E-07	&2.17696 	&1.42591E-04	&2.01696 	&1.657236\\
	                  & &1/25 	&1.98760E-07	&1.92321 	&5.64742E-05	&2.07534 	&5.345355\\	
 \bottomrule
 \end{tabular}
 \end{table}

\begin{table}[h]
\centering
 \caption{$\theta=0.1$, $\Delta t=1/200$, $H=\sqrt{2}\hat{H}$ and $h=\sqrt{2}\hat{h}$ }\label{tab6}
 \begin{tabular}{*{9}{c}}
 \midrule
   %\bottomrule
 %  &&&&&&\\
$\gamma$& $\hat{H}$&$\hat{h}$&	$\|u-U_h\|$	&Order&$\|\sigma-\Sigma_h\|$& Order&	$\frac{\text{CPU time}}{\text{Seconds}}$\\
\midrule
0.1 &     1/3  	   &1/9  	&2.06411E-06&&		4.70651E-04	&&	0.366366\\
  & 1/4  	  & 1/16 	&6.46351E-07	&2.01805 	&1.49109E-04	&1.99777 	&0.786107\\
  & 1/5  	   &1/25 	&2.56619E-07	&2.06986 	&5.89125E-05	&2.08077 	&2.066031\\
						
   20 &1/3&  	   1/9 & 	1.64081E-06&&		4.55079E-04	&&	0.349122\\
   &1/4  	   &1/16 	&4.68907E-07	&2.17696 	&1.42591E-04	&2.01696 	&0.788809\\
   &1/5  	   &1/25 	&1.98760E-07&	1.92321 	&5.64742E-05&	2.07534 	&2.07644\\
		\hline
MFE method&&&$\|u-u_h\|$&&$\|\sigma-\sigma_h\|$&&&	\\
\hline
	0.1   & &1/9  	&2.06411E-06	&&	4.70651E-04	&&	0.467958\\
	   &&1/16 	&6.46351E-07	&2.01805 	&1.49109E-04&	1.99777 	&1.692304\\
	   &&1/25 	&2.56619E-07	&2.06986 	&5.89125E-05&	2.08077 	&5.398899\\						
	  20  &&1/9  	&1.64081E-06&&		4.55079E-04	&&	0.474875\\
	  && 1/16 	&4.68907E-07&	2.17696 	&1.42591E-04&	2.01696 &	1.672167\\
	   &&1/25 	&1.98760E-07	&1.92321 	&5.64742E-05&	2.07534 &	5.473512\\	
 \bottomrule
 \end{tabular}
 \end{table}

\begin{table}[h]
\centering
 \caption{$\theta=0.3$, $\Delta t=1/200$, $H=\sqrt{2}\hat{H}$ and $h=\sqrt{2}\hat{h}$}\label{tab7}
 \begin{tabular}{*{9}{c}}
 \midrule
   %\bottomrule
 %  &&&&&&\\
$\gamma$& $\hat{H}$&$\hat{h}$&	$\|u-U_h\|$	&Order&$\|\sigma-\Sigma_h\|$& Order&	$\frac{\text{CPU time}}{\text{Seconds}}$\\
\midrule
0.1 &  1/3  &	   1/9  &	2.06411E-06&&		4.70651E-04&&		0.343843\\
  & 1/4  	  & 1/16 	&6.46351E-07	&2.01805 	&1.49109E-04	&1.99777 	&0.796913\\
   &1/5  	   &1/25 	&2.56619E-07	&2.06986 	&5.89125E-05	&2.08077 	&2.091595\\
					
  20 & 1/3  	   &1/9  	&1.64081E-06		&&4.55079E-04		&&0.344435\\
   &1/4  	   &1/16 	&4.68907E-07	&2.17696 	&1.42591E-04	&2.01696 	&0.783419\\
   &1/5  	  & 1/25 	&1.98760E-07	&1.92321 	&5.64742E-05	&2.07534 	&2.103889\\
	\hline
MFE method&&&$\|u-u_h\|$&&$\|\sigma-\sigma_h\|$&&&	\\
\hline							
	0.1     &&1/9 & 	2.06411E-06&&		4.70651E-04	&&	0.477572\\
	  && 1/16 	&6.46351E-07	&2.01805 	&1.49109E-04	&1.99777 	&1.650025\\
	   &&1/25 	&2.56619E-07&	2.06986 	&5.89125E-05	&2.08077 	&5.515452\\						
	20    &&1/9  	&1.64081E-06	&&	4.55079E-04	&&	0.481473\\
	  && 1/16 	&4.68907E-07	&2.17696 	&1.42591E-04&	2.01696 &	1.642476\\
	  && 1/25 	&1.98760E-07	&1.92321 	&5.64742E-05	&2.07534 	&5.512835\\
 \bottomrule
 \end{tabular}
 \end{table}

\begin{table}[h]
\centering
 \caption{$\theta=0.5$, $\Delta t=1/200$, $H=\sqrt{2}\hat{H}$ and $h=\sqrt{2}\hat{h}$}\label{tab8}
 \begin{tabular}{*{9}{c}}
 \midrule
   %\bottomrule
 %  &&&&&&\\
$\gamma$& $\hat{H}$&$\hat{h}$&	$\|u-U_h\|$	&Order&$\|\sigma-\Sigma_h\|$& Order&	$\frac{\text{CPU time}}{\text{Seconds}}$\\
\midrule
0.1 &     1/3  &	   1/9  &	2.06411E-06&&		4.70651E-04&&		0.339484\\
  & 1/4  	  & 1/16 	&6.46351E-07	&2.01805 	&1.49109E-04	&1.99777 	&0.790358\\
   &1/5  &	   1/25 &	2.56619E-07	&2.06986 	&5.89125E-05	&2.08077 	&2.095047\\
					
   20&1/3  &	   1/9  &	1.95560E-06	&&	4.73633E-04	&&	0.348791\\
  & 1/4  	  & 1/16 	&6.06713E-07	&2.03419 	&1.49540E-04	&2.00372 &	0.797906\\
  & 1/5  	  & 1/25 	&2.39347E-07	&2.08417 	&5.90755E-05&	2.08106 	&2.136019\\
	\hline
MFE method&&&$\|u-u_h\|$&&$\|\sigma-\sigma_h\|$&&&	\\
\hline	
						
	  0.1   & &1/9  	&2.06411E-06	&&	4.70651E-04		&&0.471017\\
	   &&1/16 	&6.46351E-07	&2.01805 	&1.49109E-04	&1.99777 	&1.694359\\
	   &&1/25 	&2.56619E-07	&2.06986 	&5.89125E-05	&2.08077 	&5.465744\\						
	20   &&1/9  	&1.95560E-06		&&4.73633E-04	&&	0.480821\\
	   &&1/16 	&6.06713E-07	&2.03419 &	1.49540E-04&	2.00372 	&1.684911\\
	   &&1/25 	&2.39347E-07&	2.08417 	&5.90755E-05&	2.08106 &	5.520902\\
 \bottomrule
 \end{tabular}
 \end{table}

\textbf{Example 4.3}
\\
Here we consider the third example to test the temporal convergence rate of second-order $\theta$ schemes. Now we choose
in (\ref{0.1}) the nonlinear term $f(u)=u^3-u$, the source term $g(x,t)=(4\pi^2-2+16\gamma\pi^4)\exp(-t)\sin(2\pi x)
 +\exp(-3t)\sin^3(2\pi x)$ and the exact solution $u(x,t)=\exp(-t)\sin(2\pi x)$ based on the taken spatial domain $\bar{\Omega}=[-1,1]$ and the temporal interval $[0,1]$. In Table \ref{tab9}, we list the error results and convergence order with CPU time, which show our TGMFE method can achieve the same calculation accuracy with approximating time second-order convergence rate as the ones by classical nonlinear Galerkin MFE method, and reduce the computing time.
\begin{table}[h]
\centering
 \caption{Temporal convergence rate}\label{tab9}
 \begin{tabular}{*{10}{c}}
 \midrule
   %\bottomrule
 %  &&&&&&\\
$\theta$&$\gamma$&$\Delta t$& $H$&$h$&	$\|u-U_h\|$	&Order&$\|\sigma-\Sigma_h\|$& Order&	$\frac{\text{CPU time}}{\text{Seconds}}$\\
\midrule	
0.1&	1&	   1/5 &   	   1/70 &	   1/4900 &  	1.17772E-05	    &           &4.68596E-04	&	        &0.3236\\
   &     &     1/10&   		    &             &     3.18843E-06  	&1.88507 	&1.29492E-04	&1.85548 	&0.4184\\
   &     &     1/20&   		    &	          &     5.91163E-07     &2.43122    &2.70235E-05	&2.26058    &0.6056\\
0.3&  10&      1/5 &  	   1/120&	   1/14400&	    1.20679E-06	    &           &4.80525E-05	&           &0.9459\\
   &    &      1/10&	        &		      &     3.24637E-07     &1.89427 	&1.32293E-05    &1.86087 	&1.2179\\
   &    &      1/20& 	        &		      &     8.07215E-08   	&2.00780 	&2.79307E-06	&2.24382 	&1.7510\\									

\hline
MFE& method     &&&&$\|u-u_h\|$&&$\|\sigma-\sigma_h\|$&&	\\
\hline								
 0.1&   1 &     1/5 &            &	   1/4900 &     1.17746E-05		&           &4.68493E-04	&	        &0.5560\\
   &     &     1/10&            &		      &     3.18549E-06	    &1.88608 	&1.29376E-04    &1.85646    &0.8439\\
   &     &     1/20&            &		      &     5.88125E-07	    &2.43732 	&2.69034E-05	&2.26571 	&1.4661\\
									
0.3&   10   &1/5       &            &      1/14400&  	1.20676E-06	    &           &4.80512E-05	&           &1.6176\\
   &    &1/10      &            &	          &     3.24601E-07     &1.89440 	&1.32278E-05	&1.86100    &2.5707\\
   &    &1/20      &            &	          &     8.06300E-08     &2.00928 	&2.78874E-06	&2.24589 	&4.4280\\
\bottomrule
 \end{tabular}
 \end{table}									
In summary, we compute three numerical examples including two-dimensional problem with exact solution, two-dimensional example covering initial value and a one-dimensional equation to test our numerical methods. In view of these calculated data in Tables \ref{tab1}-\ref{tab9}, ones know that the $\theta$ scheme can keep the second-order approximate accuracy in time, TGMFE scheme reduces the CPU time greatly without loss of calculation accuracy.

\section{Conclusions}
In this paper, we mainly consider TGMFE method and standard nonlinear Galerkin MFE method with some second-order $\theta$-schemes to solve nonlinear fourth-order reaction-diffusion equation, discuss the numerical theories including stability and error estimates, and make some numerical calculations. By the comparisons between two numerical methods, we find that TGMFE method can save the CPU-time, also see that time convergence rate is $2$ and second-order $\theta$-schemes include general Crank-Nicolson scheme and second-order backward difference method.
\par
In another work, we will apply the second-order $\theta$ scheme with finite element method to solving the following time-dependent nonlinear thermistor system
\begin{equation}\begin{split}\label{1.111111}
 u_{t}-\bigtriangleup u=&\sigma(u)|\nabla \phi|^2, \\
-\nabla\cdot(\sigma(u)\nabla \phi)=&0.
\end{split}\end{equation}
\section*{Acknowledgements} This work is supported by National Natural Science Fund (11661058, 11761053), Natural Science Fund of Inner Mongolia Autonomous Region (2016MS0102, 2017MS0107), Program for Young Talents of Science and Technology in Universities of Inner Mongolia Autonomous Region (NJYT-17-A07).

\end{document}